\RequirePackage{fix-cm}
\documentclass[smallextended]{svjour3}       
\journalname{}

\usepackage{amsmath,mathtools,amsfonts}
\usepackage{graphicx}
\usepackage{relsize,exscale}

\usepackage{tikz}
\usetikzlibrary{calc}

\usepackage[plainpages=false, colorlinks=true, citecolor=blue, 
    filecolor=black, linkcolor=blue, urlcolor=blue]{hyperref}
\usepackage[textsize=footnotesize]{todonotes}
\usepackage{bm}

\renewcommand\d[1]{\:\mathrm{d}#1}

\def\x{\mathbf{x}}
\def\y{\mathbf{y}}
\DeclarePairedDelimiter\abs{\lvert}{\rvert}%
\DeclarePairedDelimiter\norm{\lVert}{\rVert}%

\usepackage{algorithmicx}
\usepackage{algorithm}
\usepackage[noend]{algpseudocode}
\usepackage{algpseudocode}
\algrenewcommand\alglinenumber[1]{\footnotesize #1}
\algrenewcommand\algorithmicindent{0.9em}%
\algblockdefx{ForAllp}{EndForAllp}[1] {\textbf{in parallel for} #1}%

\makeatletter
\ifthenelse{\equal{\ALG@noend}{t}}%
  {\algtext*{EndForAllp}}
  {}%
\makeatother

\algblock{Input}{EndInput}
\algnotext{EndInput}
\algblock{Output}{EndOutput}
\algnotext{EndOutput}

\begin{document}

\title{Space-Fractional Diffusion with Variable Order and Diffusivity: Discretization and Direct Solution Strategies}
\titlerunning{Space-Fractional Diffusion with Variable Order and Diffusivity}

\author{Hasnaa Alzahrani \and George Turkiyyah \and Omar Knio \and David Keyes}

\institute{H. Alzahrani \at King Abdullah University of Science and Technology, Thuwal, Saudi Arabia \and 
G. Turkiyyah \at King Abdullah University of Science and Technology, Thuwal, Saudi Arabia \and
O. Knio \at King Abdullah University of Science and Technology, Thuwal, Saudi Arabia \\ \email{omar.knio@kaust.edu.sa} \and 
D. Keyes \at King Abdullah University of Science and Technology, Thuwal, Saudi Arabia}

\date{}

\maketitle

\begin{abstract}
We consider the multidimensional space-fractional diffusion equations with spatially varying diffusivity and fractional order. Significant computational challenges are encountered when solving these equations due both to the kernel singularity in the fractional integral operator and to the resulting dense discretized operators, which quickly become prohibitively expensive to handle because of their memory and arithmetic complexities. 

In this work, we present a singularity-aware discretization scheme that regularizes the singular integrals through a singularity subtraction technique adapted to the spatial variability of diffusivity and fractional order. This regularization strategy is conveniently formulated as a sparse matrix correction that is added to the dense operator, and is applicable to different formulations of fractional diffusion equations.
We also present a block low rank representation to handle the dense matrix representations, by exploiting the  ability to approximate blocks of the resulting formally dense matrix by low rank factorizations. A Cholesky factorization solver operates directly on this representation using the low rank blocks as its atomic computational tiles, and achieves high performance on multicore hardware. 

Numerical results show that the singularity treatment is robust, substantially reduces discretization errors, and attains the first-order convergence rate allowed by the regularity of the solutions. They also show that considerable savings are obtained in storage ($O(N^{1.5})$) and computational cost ($O(N^2)$) compared to dense factorizations. This translates to orders-of-magnitude savings in memory and time on multi-dimensional problems, and shows that the proposed methods offer practical tools for tackling large nonlocal fractional diffusion simulations.  


\end{abstract}

\section{Introduction}
\label{sec:intro}

Simulations involving space-fractional diffusion operators are becoming increasingly important in a number of application domains.  Their ability to model phenomena of anomalous nonlocal diffusion in fractured and granular media, and account for long range interactions beyond classical Brownian motion, have made them powerful tools in several areas of relevance to industrial and environmental applications \cite{delia20a,lucchesi20}. Beyond anomalous diffusion, diverse problems in physics-informed neural networks \cite{karniadakis19}, image denoising \cite{gilboa09,delia21a}, and sampling from random gaussian fields in spatial statistics \cite{bolin18} also benefit from the ability of fractional operators to capture nonlocal effects and control solution regularity. In many of these applications, the heterogenous case, with spatially-varying fractional order and diffusivity, is particularly useful. As a result, there is significant interest in the development of fast and accurate methods for the solution of variable coefficient problems.

There are however two primary challenges in fractional diffusion simulations, when attempted on realistic problems in multiple spatial dimensions. The first one comes from the singularity of the kernel in the formulation of the fractional integral operator. Singularity of the kernel implies that standard quadrature rules will converge slowly and will not obtain the convergence rates that can be expected with smooth integrands. For the special case of the fractional Laplacian, i.e., problems with constant fractional order and constant diffusivity, treatments of the singularity have been proposed. A singularity subtraction method to regularize the integral is described in \cite{pozrikidis16}. This method was substantially enhanced in \cite{minden20} where the singularity subtraction is modified and limited to a local neighborhood through a radial windowing function that is shown to be quite effective computationally in 2D and 3D. \cite{darve_xu20} uses a similar method in the context of an isogeometric 2D discretization and shows linear convergence with mesh size, a rate that cannot be obtained by a finite element discretization that does not explicitly treat the singularity \cite{acosta17}. A finite difference method for the constant coefficient fractional Laplacian, which splits the kernel function into two weakly-singular parts, is introduced in \cite{duo19,duo18a} as a fractional analogue of the classical central difference schemes. A finite difference method for a variable diffusivity problem in 1D is described in \cite{mustapha20}. There is still however no general treatment for variable order and variable diffusivity in fractional multidimensional problems.

The second challenge comes from the fact that the discretization of the integral operator results in a dense matrix, a consequence of the non-local nature of fractional diffusion. If not effectively tackled, this imposes prohibitive computational requirements in both memory and runtime on the numerical simulations. For constant fractional order and diffusion coefficients, these challenges may be addressed by taking advantage of the homogeneity to render the problem tractable. For example, when using a regular uniform discretization, the construction and storage of only a small representative translation-invariant portion of the problem is sufficient, because of the block-Toeplitz with Toeplitz blocks structure of the resulting matrix \cite{du15,minden20}. Consequently, application of the discrete operator can be performed in $O(N \log N)$ using FFT. For variable fractional order in 1D, \cite{jia20} approximates the discrete operator by a scaled sum of Toeplitz matrices, which permits its application in log-linear asymptotic complexity. However, in the general multidimensional case when either the order or the diffusion coefficient is not homogeneous, or when the discretization is unstructured, or when the simulation domain is bounded, the Toeplitz structure is not helpful and alternative representations are needed, especially as the discretizations gets refined.


Hierarchical ($\mathcal{H}$) matrices provide such an alternative and have been shown to be effective and general-purpose representations for the discretizations of fractional operators in one and two dimensional problems \cite{zhao17,massei19,karkulik19,xu18,boukaram20}. Hierarchical matrix representations allow substantial and accuracy-tunable compression of the dense matrix, by approximating certain blocks of the matrix with low rank factorizations. These low rank blocks are not necessarily of the same size but can be of different granularity representing different levels of the hierarchical representation. Hierarchically low rank matrices reduce the $O(N^2)$ memory footprint of dense matrices to $O(k N \log N)$ or even $O(k N)$ in the case of the $\mathcal{H}^2$ representation with nested bases, where $k \ll N$ is a representative block rank that depends on the desired quality of the approximation.  $\mathcal{H}$- and $\mathcal{H}^2$-matrices also allow operator application to be performed with similar log-linear and linear complexity, respectively. Iterative solvers can then be readily built using the hierarchical matrix representation, as the inner kernel of an iterative solver is a matrix-vector multiplication that can be efficiently performed.

While iterative methods, when preconditioned by appropriate preconditioners, can lead to effective solvers, there are many scenarios in which direct solvers are desirable. Direct factorization-based solvers are more robust, require no parameter tuning, can effectively handle simulations with multiple right hand sides at negligible additional cost beyond the initial factorization, and can be updated via Sherman-Morrison-Woodbury formulas when local modifications in the form of low rank updates are made to the problem. 
Unfortunately, direct solvers for general hierarchical matrices are not particularly efficient. While simple blocking structures (HODLR and HSS weak-admissibility structures) allow for direct factorization methods, the rank growth in the off-diagonal blocks of these representations leads to undesirable growth in the runtime. Direct factorizations of general hierarchical matrices with strong-admissibility blocking have large constants in their complexity estimates and introduce data dependencies to make them impractical, particularly on modern multicore hardware. It is therefore desirable to use alternative representations that provide sufficient memory compression to be able to store discretizations of simulations of practical interest, and achieve high performance on direct factorizations.

In this paper, we propose strategies for addressing the two difficulties outlined above in the inhomogeneous multidimensional case. In particular, we propose a singularity treatment that allows a convenient quadrature rule, such as the trapezoidal rule, to be used in the discretization of the integral operators. Our method uses a singularity subtraction method that takes into account the spatially varying coefficients, generalizing the treatment of the constant coefficient case in \cite{minden20}. We show it to be effective for different formulations of fractional diffusion, and attains the first-order convergence allowed by the solution, which does not generally have sufficient regularity for higher order convergence. We also present a 
practical block low rank matrix representation of the resulting operator and show that it achieves substantial memory reduction of $O(k N^{1.5})$ compared to an $O(N^2)$ dense representation. This representation is the foundation of a direct solver that uses the low rank blocks as its computational tiles. The solver relies on randomized sampling to produce a block low rank Cholesky factorization in only $O(k N^2)$ operations, and with substantial concurrency in its operations. Even though the asymptotic rate is not optimal, in practice, and for many problems of relevant size, the computations are well tolerated and balanced by more efficient execution on multi/many-core architectures.

The rest of the paper is organized as follows. Section 2 describes the formulations of the fractional diffusion operators that we consider. Section 3 describes the singularity subtraction treatment for spatially varying fractional order and diffusivity in the simple one-dimensional context. Section 4 presents the singularity treatment in the multidimensional case. Section 5 describes the matrix representation that uses low rank factorizations in matrix blocks and the direct solver that operates on the compressed representation to generate a Cholesky factorization. Numerical results are presented in Section 6. They show the convergence attained by the discretization strategy in 1D and 2D with various variable diffusivity and fractional order examples, as well as the memory compression produced by the blockwise low rank representation and the runtime savings achieved by the solver on large representative 2D problems. We conclude and outline future work in Section 7.

\section{Formulation with variable diffusivity and fractional order}
\label{sec:formulation}

\begin{figure}
  \begin{center}
    \includegraphics[width=.8\textwidth]{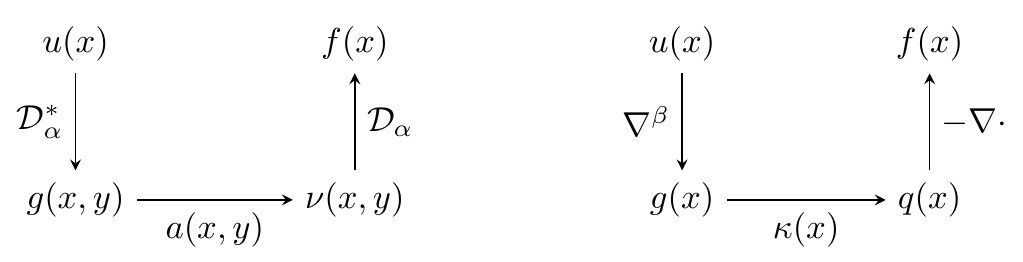}
  \end{center}
  \caption{Nonlocal diffusion frameworks. The symmetric formulation on the left uses a nonlocal divergence operator and its adjoint as a nonlocal gradient. The non-symmetric formulation on the right uses a nonlocal flux and uses the classical local divergence.}   
  \label{fig:form}
\end{figure}

Different formulations of nonlocal diffusion are possible and have been proposed in the literature. They differ in the way that the fractional gradient and/or 
divergence are defined.  In this work, we consider two specific formulations, both consistent with nonlocal conservation laws; see Fig.~\ref{fig:form} for an illustration. 

The first formulation considered follows a nonlocal mechanics formalism~\cite{DuEtAl2013,du19}, in which a nonlocal generalization of divergence is introduced and its adjoint is used as a nonlocal gradient. In the second formulation, a nonlocal flux is first defined using a fractional gradient, and the classical divergence the fractional flux is used to defined the fractional diffusion term.  In light of these definitions, we shall refer to these two formulations
as symmetric and non-symmetric, respectively.

\subsection{Symmetric formulation}

In the symmetric case, the fractional diffusion term is defined in terms of (i) the generalized (non-local) divergence operator~\cite{DuEtAl2013},
\begin{equation}
\mathcal{D}(\mathbf{\nu})(\x) \coloneqq \int_{\mathbb{R}^n} \left( \nu(\x, \y) + \nu(\y, \x) \right) \cdot \alpha(\x, \y) dy 
\end{equation}
where $\nu(\x, \y):\ \mathbb{R}^n \times \mathbb{R}^n \to \mathbb{R}^n$ is a two-point vector-field $\x$ and $\y$ are points in $\mathbb{R}^n$, and 
$\mathbf{\alpha}(\x, \y):\ \mathbb{R}^n \times \mathbb{R}^n \to \mathbb{R}^n$ is an antisymmetric vector field,  i.e.\ $\mathbf{\alpha}(\x, \y) = - \mathbf{\alpha}(\y, \x)$, and (ii) its adjoint,
\begin{equation}
  \mathcal{D}^*(u)(\x, \y) = - \left( u(\y) - u(\x) \right) \mathbf{\alpha} (\x, \y)
\end{equation}
where $u(\x)$ is a scalar field over $\mathbb{R}^n$.  
\begin{equation}
  \mathcal{D}^*(u)(\x, \y) = - \left( u(\y) - u(\x) \right) \mathbf{\alpha} (\x, \y) ,
\end{equation}
which is viewed as a generalized (nonlocal) gradient.  Using these definitions, and given a second-order tensor, $a$, satisfying $a(\x, \y) = a(\y, \x)$ and 
$a = a^T$, one defines the generalized diffusion term:
\begin{equation}
  \mathcal{D} \left( a \cdot \mathcal{D}^* \right) (\x) = -2 \int_{\mathbb{R}^n} \left( u(\y) - u(\x) \right) \gamma(\x, \y) d\y
  \label{eq:integral}
\end{equation}
where 
$$
\gamma(\x,\y) \coloneqq \mathbf{\alpha}(\x, \y) \cdot a(\x, \y) \cdot \mathbf{\alpha}(\x, \y) .$$

In this work, we shall suppose that $a$ is diagonal, i.e., it reduces to a scalar of the form
\begin{equation} 
a(\x, \y) = \sqrt{\kappa(\x) \kappa(\y)} ,
\label{eq:forma}
\end{equation}
with $\kappa(\x) \ge \delta > 0$, $\delta$ a constant independent of $\x$, and that $\mathbf{\alpha}$ is given by:
\begin{equation}
\mathbf{\alpha}(\x, \y) \equiv \frac{\y-\x}{| \y- \x |^{\frac{n}{2} + \frac{\beta(\x)+\beta(\y)}{2} + 1}} ,
\end{equation}
with $0 < \beta(\x) < 1$.  This yields 
\begin{align}
\gamma(\x, \y) &= \frac{\sqrt{\kappa(\x)\kappa(\y)} \ |\y-\x|^2}{|\y-\x|^{d+\beta(\x)+\beta(\y)+2}} 
             = \frac{\sqrt{\kappa(\x)\kappa(\y)}}{|\y-\x|^{n+\beta(\x)+\beta(\y)}} ,
\label{eq:kernel}
\end{align}
and enables us to interpret
\begin{equation}
  \mathcal{L}_\beta u(\x) \coloneqq -2 \int_{\mathbb{R}^n} \left( u( \y ) - u( \x ) \right) \frac{\sqrt{\kappa(\x)\kappa(\y)}}{|\y-\x|^{n+\beta(\x)+\beta(\y)}} d\y
  \label{eq:Lsym}
\end{equation}
as a fractional (sub-) diffusion operator of variable order $\beta( \x)$, and variable diffusivity field $\kappa(\x)$.  Note that the expression
of $a$ given in~(\ref{eq:forma}) is appealing on dimensional grounds, as $\kappa(\x)$ has dimension $[{\rm L}]^{2\beta(\x)} / {\rm T}$ where
$[{\rm L}]$ and $[{\rm T}]$ respectively denote length and time dimensions.

\begin{figure}[h]
  \begin{center}
    \includegraphics[width=.35\textwidth]{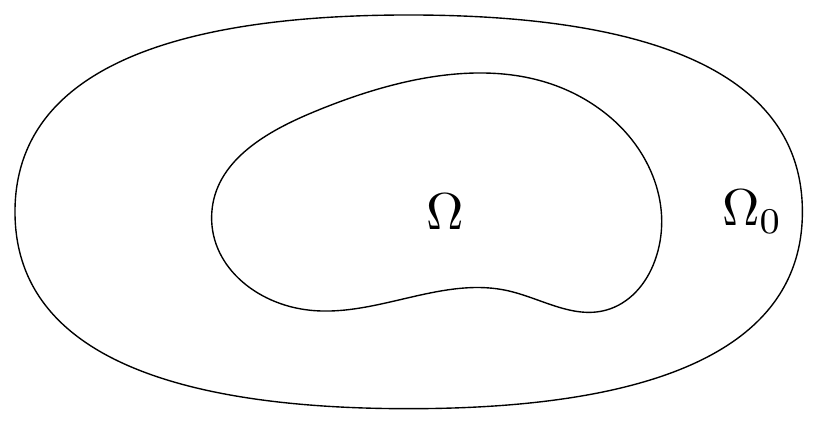}
  \end{center}
  \caption{Schematic illustration of the solution domain, $\Omega$, and the surrounding region $\Omega_0$ over which volume constraints are imposed.} 
  \label{fig:sketch}
\end{figure}

We shall generally focus on solving the fractional diffusion equation,
\begin{equation}
  \mathcal{L}_\beta u (\x) = f(\x),
  \label{eq:fheatsym}
\end{equation}
in a domain $\Omega \subset \mathbb{R}^n$, where $f$ is a given source term.  In case there are no interactions with the region outside the domain, 
the integral in (\ref{eq:Lsym}) is restricted to $\Omega$.  In a more general setting, we consider that the domain $\Omega$ is surrounded by 
an enclosing region $\Omega_0 \subset \mathbf{R}^n$ with $\Omega \cap \Omega_0 =  \emptyset$, and that nonlocal interactions occur between 
points in $\Omega$ and $\Omega_0$; see Fig.~\ref{fig:sketch}.  In this case, the integral (\ref{eq:Lsym}) is performed over $\Omega \cup \Omega_0$, and the fields $u(\y)$, 
$\kappa(\y)$ and $\beta(\y)$ are assumed to be specified for $\y \in \Omega_0$.  These ``volume constraints'' are the nonlocal generalizations of the classical (local) Dirichlet or Neumann boundary conditions \cite{du12,delia20b}.  In the following, we shall focus on homogeneous Dirichlet type conditions, with $u(\y) = 0$
when $\y \in \Omega_0$.

We also note that the standard fractional Laplacian may be recovered as a special case of the above, namely for constant $\beta$ and $\kappa$, and $\Omega \cup \Omega_0 = \mathbb{R}^n$.

\subsection{Non-symmetric formulation}

The non-symmetric formulation consists in first defining a nonlocal flux, $\mathbf{Q}^\beta$, based on a fractional derivative of order 
$\beta$~\cite{samko13,pozrikidis16,kwasnicki17,lischke20,lucchesi20,boukaram20}.  Allowing for variable order and diffusivity, we set
\begin{equation}
\mathbf{Q}^\beta= -\kappa(\x) \mathbf{\nabla}^\beta u(\x)
\label{eq:fracflux}
\end{equation}
where $\kappa$ is the diffusivity,
\begin{equation}
\mathbf{\nabla}^\beta u(\x)= \omega(\x) \int_{\mathbb{R}^n} \frac{\y - \x}{| \y - \x | ^{\beta (\x)+n+1}} u(\y)d\y
\label{eq:fracgrad}
\end{equation}
is the fractional gradient of order $\beta$, whereas
\begin{equation}
\omega (\x):=  \frac{2^{\beta(\x)} \Gamma(\frac{n+\beta (\x)+1}{2})}{\pi^{\frac{n}{2}} \Gamma(\frac{1-\beta(\x)}{2})}
\end{equation}
is scaling factor that depends on the (variable) order and on the number of spatial dimensions.  Note that $\omega$
may be dropped by suitably rescaling the diffusivity, $\kappa$, and that with $\kappa (\x)$ having dimension 
$[L]^{\beta(\x)+1}/[T]$ the fractional flux $\mathbf{Q}^\beta$ is dimensionally homogeneous across space.

In the non-symmetric case, we shall consider the solution of the fractional diffusion equation,
\begin{equation}
{\cal N}_\beta u ( \x ) = f ( \x)
\end{equation}
where 
\begin{equation}
{\cal N}_\beta u ( \x ) \equiv - \mathbf{\nabla} \cdot \mathbf{Q}^\beta ,
\end{equation}
and $\mathbf{\nabla}$ the classical gradient of order 1.  As in the symmetric case, we consider the Dirichlet problem
in a bounded domain $\Omega$.  As in the symmetric case, this requires specifying values of $u$ outside $\Omega$, but
unlike the symmetric case the fields $\kappa$ and $\beta$ need not be specified outside the domain.  Note that in the 
special case where $u$ vanishes identically outside the domain, the integral in (\ref{eq:fracgrad}) may simply be restricted to
$\Omega$.


\section{Singularity-aware discretization in the 1D case}

We first describe our scheme for the treatment of the kernel singularity in the one-dimensional context, which we extend to the multidimensional case in the following section. We treat the variable diffusion coefficient first, followed by the variable fractional order. 

\subsection{Spatially varying nonlocal diffusion coefficient} 
\label{sec:geomean1d}

The two-parameter diffusivity coefficient $a(x, y)$ that appears in the symmetric formulation of fractional diffusion may take different algebraic forms, with a common one being the geometric mean of a one-parameter diffusivity field $\kappa$
\begin{equation}
a(x,y) = \kappa(x)^{1/2} \kappa(y)^{1/2} = c(x) c(y)
\end{equation}
where we define $c(x) := \kappa(x)^{1/2}$. The kernel function becomes
\begin{equation}
  \gamma(x, y) = \frac{c(x) c(y)}{|y-x|^{1+2\beta}}
\end{equation}
and the integral in (\ref{eq:integral}) is written as:
\begin{equation}
 \mathcal{L}[u(x)] = -2 \int_{\Omega \cup \Omega_0} \frac{[u(x) - u(y)] \, [c(x) c(y)] }{|y-x|^{1+2\beta}}   \, dy
 \label{eq:int2}
\end{equation}

\begin{figure}[ht]
  \begin{center}
    \includegraphics[width=0.4\textwidth]{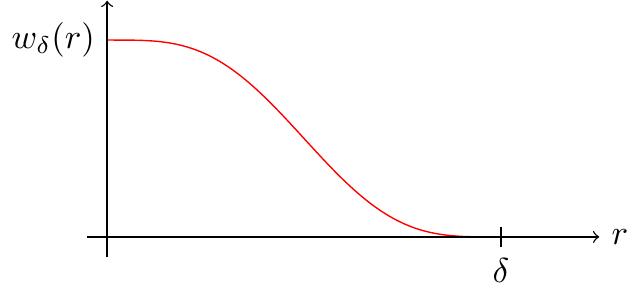}
    \caption{Windowing function $w_\delta(r) = 1 - 35(r/\delta)^4 + 84(r/\delta)^5 - 70(r/\delta)^6 + 20(r/\delta)^7$ for $r < \delta$}
    \label{fig:window}
  \end{center}
\end{figure}

It is possible to remove the singularity of the integral of (\ref{eq:int2}) at $y = x$ by adding a term to the integrand and subtracting it in a separate term that can be handled more readily. We can use the Taylor series expansions
\begin{align}
 u(y) &= u(x) + u'(x) (y-x) + \tfrac{1}{2} u''(x)(y-x)^2 + O(|y-x|^3 ) \label{eq:u_taylor} \\
 c(y) &= c(x) + c'(x) (y-x) + O(|y-x|^2)  \label{eq:c_taylor}
\end{align}
to define the following desingularization term in a small window around the point $x$:
\begin{equation}
  p_x(y) = 2 \, \underbrace{w(\abs{y-x})}_{\text{local window}} \ 
    \underbrace{\left[ u'(x) (y-x) + \tfrac{1}{2} u''(x)(y-x)^2 \right]}_{\text{local approximation of } u(y)-u(x)} 
    \underbrace{\left[ c(x)^2 + c(x) c'(x) (y-x)\right]}_{\text{local approximation of } a(x,y)=c(x)c(y)} 
  \label{eq:desing}
\end{equation}
where $w(\abs{y-x})$ is a suitably chosen radial regularization function such $w(\abs{y-x}) = 1 + O(|y-x|^4)$ as $y \to x$. A local polynomial windowing function \cite{minden20} that satisfies this condition is shown in Fig.~\ref{fig:window}. The $O(|y-x|^3)$ terms in (\ref{eq:u_taylor}) and the $O(|y-x|^2)$ terms in (\ref{eq:c_taylor}) have been dropped in (\ref{eq:desing}), either because they result in odd-power terms that integrate to zero or in higher order terms $O(|y-x|^4)$ that only introduce smooth terms to the integrand, and do not affect the discretization in what follows.

The regularized version of the integral of (\ref{eq:int2}) can then be written as: 
\begin{align}
&  \int_{\Omega \cup \Omega_0} 
\left[ \frac{-2\left[u(y) - u(x)\right] \, a(x,y)}{|y-x|^{1+2\beta}} + 
\frac{p_x(y)}{|y-x|^{1+2\beta}} \right] dy 
\, - 
\int_{\Omega \cup \Omega_0} 
\frac{p_x(y)}  {|y-x|^{1+2\beta}} \, dy
\label{eq:regularized}
\end{align}
In this form, the integrand of the first term goes to zero as $y \to x$. It is also continuously differentiable with an integrable second derivative, allowing a trapezoidal rule to be used in the discretization of the integral. 

We consider the discretization of the first integral at a point $x = x_i$ on a regular grid with spacing $h$.
The numerator $p_x(y)$ of the singularity-removing term in the first integrand, when expanded, has two terms with odd-powers of $(y-x_i)$ whose integrals vanish, and remaining terms that can be expressed in the form:
\begin{equation}
  w(y-x_i) \ c_i[2 ( c_i u'_i )' - c_i u''_i]  \ (y-x_i)^2
  \label{eq:reg_kappa1d}
\end{equation}
where $c_i = c(x_i)$, $u_i = u(x_i)$, $u'_i = u'(x_i)$, and $u''_i = u''(x_i)$. 
We define the quantity 
\begin{equation}
\overline{u''_i} \coloneqq c_i[2 (c_iu_i')'-c_iu_i''],
\end{equation}
which admits a symmetric discretization as:  
\begin{align}
\overline{u''_i}  & \approx c_i \tfrac{1}{h^2} \left[ 2(c_{i-1/2} u_{i-1} -(c_{i-1/2} + c_{i+1/2}) u_i + c_{i+1/2} u_{i+1}) - (c_i u_{i-1} - 2 c_i u_i + c_i u_{i+1}) \right] \\
&= c_i \tfrac{1}{h^2} \left[ (2c_{i-1/2} - c_i)u_{i-1} - 2 (c_{i-1/2} + c_{i+1/2} - c_i)u_i + (2c_{i+1/2} - c_i) u_{i+1} \right] 
\end{align}
and that we write as
\begin{equation}
\overline{u''_i} = \frac{1}{h^2} \left( k_{i-1} u_{i-1} - (k_{i-1} + k_{i+1}) u_i + k_{i+1} u_{i+1} \right)
\label{eq:ubar1d}
\end{equation}
with $k_{i\pm 1} \coloneqq c_i (2c_{i\pm 1/2} - c_i)$. Using $O(h^2)$ linear approximations of $c_{i\pm 1/2}$, we get that $k_{i+1} = c_i c_{i+1}$ and $k_{i-1} = c_i c_{i-1}$. We note that for the constant coefficient case, with $c_{i-1/2} = c_{i+1/2} = c_i = \kappa^{1/2}$, the expression of $\overline{u''_i}$ reduces to the usual central difference formula:
\begin{equation}
  u''_i \approx  \frac{1}{h^2} (\kappa u_{i-1} - 2 \kappa u_i + \kappa u_{i+1} )
\end{equation}
which appears when regularizing the constant coefficient fractional Laplacian \cite{minden20}. In fact, the notation $\overline{u''_i}$ was chosen to reflect the generalization to the variable diffusivity case. 

Using the trapezoidal rule, the first integral in (\ref{eq:regularized}) can be discretized at point $x_i$ as follows
\begin{equation}
 h \sum_{j\ne i} \left[ -2(u_j - u_i) \gamma_{ij} + \overline{u''_i} \frac{w(\abs{x_j - x_i})}{|x_j-x_i|^{2\beta-1}} \right]
\label{eq:I1}
\end{equation}
where $\gamma_{ij} = \gamma(x_i, x_j)$. The second integral in (\ref{eq:regularized}) can be computed to high precision by a separate adaptive quadrature since its integrand is non-zero only in a small region around $x_i$:
\begin{equation}
  - \overline{u''_i} \int_{\Omega \cup \Omega_0} \frac{w(\abs{y - x_i})}{|y-x_i|^{2\beta-1}} dy
\label{eq:I2}
\end{equation}

Assuming extended Dirichlet conditions in the region $\Omega_0$, i.e., $u_j = 0$ for grid points $j$ in $\Omega_0$, and a grid of size $N$ in the interior region $\Omega$, the final discretization of (\ref{eq:regularized}) can be written as
\begin{equation}
  (B + D + C) u
\end{equation}
where $B$ is an $N \times N$ matrix with entries
\begin{equation}
B_{ij} = 
\begin{cases}
 -2 h \gamma_{ij},  & \quad i \ne j \\
  0,  & \quad i = j
\end{cases} 
\end{equation}
$D$ is a diagonal matrix with entries 
\begin{equation}
   D_i = 2 h \sum_{\substack{j \ne i \\ j \in \Omega \cup \Omega_0}} \gamma_{ij}, 
\end{equation}
and the last term in (\ref{eq:I1}) together with (\ref{eq:I2}) contribute a tridiagonal matrix $C$ to the discretization. The tridiagonal sparsity pattern of $C$ is a consequence of (\ref{eq:ubar1d}), and is obviously the same sparsity pattern of a 3-point discretization of the classical (non-fractional) 1D Laplacian. $C$ is also symmetric leading to a symmetric discrete operator $A= B + D + C$. 

\subsection{Spatially varying fractional order} 

We now consider the case of variable order with $\beta(x)$ a function of the spatial variable. We assume $\beta(x)$ varies smoothly and has a bounded derivative $\beta'(x)$ everywhere, and $\kappa(x) = 1$.  The spatially varying kernel is then:
\begin{equation}
  \gamma(x, y) = \frac{1}{|y-x|^{1+\beta(x)+\beta(y)}}   
  \label{eq:varkernel1d}
\end{equation}
Expressing $\beta(y)$ as a Taylor series at point $x$,
\begin{equation}
  \beta(y) = \beta(x) + \beta'(x) (y-x) + O(|y-x|^2)
  \label{eq:betataylor}
\end{equation}
allows us to express the kernel (\ref{eq:varkernel1d}) as the expansion:
\begin{equation}
   \gamma(x, y) = \frac{1}{|y-x|^{1+\beta(x)+\beta(y)}} = \frac{1}{|y-x|^{1+2\beta(x)}} (1 - \beta'(x) (y-x) 
  \log |y-x| + O(|y-x|^2))
  \label{eq:kernelseries}
\end{equation}
Combining (\ref{eq:kernelseries}) with the Taylor series expansion
\begin{equation}
  u(y) = u(x) + u'(x) (y-x) + \tfrac{1}{2}u''(x) (y-x)^2 + O(|y-x|^3)
  \label{eq:taylor1d}
\end{equation}
allows us to subtract, locally around $y=x$, the variable order singularity in the original (\ref{eq:integral}) and write the integral as:
\begin{subequations}
\label{eq:varorder}
\begin{align}
 & \int \left[  \frac{-2(u(y) - u(x))}{|y-x|^{1+\beta(x)+\beta(y)}} + \frac{w(\abs{y-x}) \left(  u''(x) (y-x)^2 - 2 u'(x) \beta'(x) (y-x)^2 \log|y-x| \right)}{|y-x|^{1+2\beta(x)}} \right] \ dy \label{eq:smooth} \\[3pt]
 & -u''(x) \int \frac{w(\abs{y-x})}{|y-x|^{2\beta(x)-1}} dy \ + \ 2u'(x) \int  \frac{w(\abs{y-x}) \beta'(x) \log|y-x|}{|y-x|^{2\beta(x)-1}}dy
 \label{eq:analytic}
\end{align}
\end{subequations}
The $O(|y-x|^2)$ terms in (\ref{eq:kernelseries}) and the $O(|y-x|^3)$ terms in (\ref{eq:taylor1d})  have been dropped in (\ref{eq:varorder}), either because they result in odd-power terms that integrate to zero or in higher order terms $O(|y-x|^4)$ that only introduce smooth terms to the integrand that do not affect the discretization of the regularized integral. 

As with the variable diffusivity case, the regularizing term of (\ref{eq:smooth}) was chosen to make the integrand go to zero as $y \to x$ with an integrable second derivative, allowing a trapezoidal rule to be used to evaluate the integral. The two integrals of (\ref{eq:analytic}) can be computed either analytically when the window function $w$ has a simple form or by an adaptive quadrature method.

At a point $x_i$, the discretization of (\ref{eq:smooth}) is then of the form: 
\begin{equation}
  h \sum_{j\ne i} \left[
   \frac{-2(u_j - u_i)}{|x_j - x_i|^{1 + \beta_i+\beta_j}}
   + u''_i \frac{w(\abs{x_j - x_i})}{|x_j - x_i|^{2\beta_i-1}}
   - 2u'_i \frac{w(\abs{x_j - x_i}) \beta'_i \log |x_j - x_i|}{|x_j - x_i|^{2\beta_i-1}}
  \right]
\end{equation} 
where  $\beta_i = \beta(x_i)$, $\beta_j = \beta(x_j)$, and $\beta'_i = \beta'(x_i)$. $u'_i$ and $u''_i$ are the first and second derivatives at $x_i$ and may be approximated by the usual $O(h^2)$ finite difference formulas $u''_i \approx \tfrac{1}{h^2}(u_{i+1} - 2u_i + u_{i-1})$ and $u'_i \approx \tfrac{1}{2h}(u_{i+1} - u_{i-1})$. The discretization of the rest of (\ref{eq:varorder}) is of the form: 
\begin{equation}
  - u''_i I_{ai} +2 u'_i I_{bi}
\end{equation}
where $I_{ai}$ and $I_{bi}$ are the two integrals of (\ref{eq:analytic}) which can be evaluated separately at every point $x = x_i$ with an appropriate quadrature. 

As with the variable diffusivity case, the final discretized operator is also of the form $A = B + D + C$ where $B$ is a dense matrix whose entries involve kernel evaluations, $D$ is a diagonal matrix that includes extended Dirichlet conditions in $\Omega_0$, and $C$ is a tridiagonal matrix resulting from the regularization of the integral. $C$ is not formally symmetric, however it can be replaced to $O(h^2)$ accuracy by its symmetrized version $(C + C^T)/2$.

The non-symmetric integral formulation of fractional diffusion described in Section~\ref{sec:formulation} may also be regularized using a similar strategy. The details are described in Appendix \ref{app:nonsym}.


\section{Singularity-aware discretization in the multidimensional case}
\label{sec:singularity_nd}

The regularization strategy in the multi-dimensional case is conceptually similar to the one-dimensional setting. The fractional operator is now  
\begin{equation}
   \mathcal{L}[u(\x)] = -2 \int_{\Omega \cup \Omega_0} \frac{u(\y) - u(\x)}{|\y - \x|^{n + \beta(x) + \beta(y)}} a(\x, \y) \d \y
  \label{eq:integral_nd}
\end{equation}

In order to discretize the integral in (\ref{eq:integral_nd}), we subtract the singularity so as to obtain a sufficiently regular integrand that allows an $n$-dimensional trapezoidal rule to be used, and handle the singularity-correction term by a separate quadrature. We first consider the spatially varying coefficient $\kappa(\x)$ and then the spatially varying fractional order $\beta(\x)$. 

\subsection{Variable diffusion coefficient}
\label{sec:varkappa_nd}

We consider the case where the two-argument nonlocal diffusion coefficient $a(\x, \y)$ is defined as the geometric mean of a diffusion coefficient $\kappa(x)$, so we write $a(\x, \y) = \kappa(\x)^{1/2} \, \kappa(\y)^{1/2}$. The case where the function $a(\x, \y)$, which must be symmetric, is the arithmetic mean of its two constituents, or has other forms, can be handled in a similar fashion, and we skip the details. 

We consider the evaluation at a point $\x = \x_i$ and make use of the following Taylor series expansions for $u(\y)$ and $c(\y) := \kappa(\y)^{1/2}$, around $\x_i$:
\begin{align}
  u(\y) & = u(\x_i) + \nabla u(\x_i)^T (\y - \x_i) + \tfrac{1}{2} (\y-\x_i)^T \nabla^2 u(\x_i) (\y-\x_i) + \cdots \label{eq:taylor_nd} \\
  c(\y) & = c(\x_i) + \nabla c(\x_i)^T (\y - \x_i) + \cdots
\end{align}
The term needed to cancel the singularity of (\ref{eq:integral_nd}) at $\y = \x_i$ takes the form
\begin{equation}
  C(\y) = \frac{2 \, w(\norm{\y-\x_i}) \ \bm{[} 
   \nabla u_i^T (\y - \x_i) + \tfrac{1}{2} (\y-\x_i)^T \nabla^2 u_i (\y-\x_i) \bm{]}  \ 
    [c_i^2 + c_i \nabla c_i^T (\y - \x_i) ] }{\norm{\y-\x_i}^{n+2\beta}}
\end{equation} 
which may be simplified, after removing terms involving odd powers of $(\y-\x_i)$ whose integrals vanish, to:
\begin{equation}
  C(\y) = \frac{ w(\norm{\y-\x_i}) \ (\y-\x_i)^T \bm{[} 2 c_i \nabla \cdot (c_i \nabla u_i) - c_i^2 \nabla^2 u_i \bm{]} (\y-\x_i) } 
  {\norm{\y-\x_i}^{n+2\beta}}
  \label{eq:reg_kappadd}
\end{equation}
which involve both the Laplacian and a variable coefficient Laplacian of $u(\x)$ at $\x_i$. Simplifying further, by noting that terms that involve cross products of the different components of the $n$-dimensional vector $(\y - \x_i)$ also have integrals that vanish, we can write the numerator of (\ref{eq:reg_kappadd})  as a sum of  terms involving the derivatives of $u(\x)$ in the $n$ coordinate directions:
\begin{equation}
  \sum_{d=1}^n w(\norm{\y-\x_i}) \left[ 2 c_i \partial_d (c_i \partial_d u_i) -c_i^2 \partial_{dd}^2 u_i \right] (\y_d - \x_{i,d})^2
  \label{eq:reg_kappa2d_num}
\end{equation}
where $\y_d$ and $\x_{i,d}$ are the $d$-th components of the $\mathbf{R}^n$ vectors $\y$ and $\x_i$, respectively, and $\partial_d$ is the derivative in the $d$-th direction. Each of the $n$ summands in (\ref{eq:reg_kappa2d_num}) looks like the one-dimensional singularity removing term in (\ref{eq:reg_kappa1d}), and we can therefore use a similar discretization to the one in Section \ref{sec:geomean1d}, to write $C(\y)$ as a sum of $n$ terms, each corresponding to a coordinate direction:
\begin{equation}
  C(\y) = \sum_{d=1}^n  \overline{\partial^2_{dd} u_i} \, \frac{ w(\norm{\y-\x_i}) (\y_d - \x_{i,d})^2 } {\norm{\y-\x_i}^{n+2\beta}}
\end{equation}
where  $\overline{\partial^2_{dd} u_i} \coloneqq 2 c_i \partial_d (c_i \partial_d u_i) -c_i^2 \partial_{dd}^2 u_i$, and can be conveniently discretized on a regular grid similarly to (\ref{eq:ubar1d}), interpreted in the $d$-th coordinate, i.e., with the $i+1$ and $i-1$ subscripts referring to the next and previous grid points in the coordinate direction $d$, respectively:
\begin{equation}
  \overline{\partial^2_{dd} u_i} \approx  \tfrac{1}{h^2} \left( k_{i-1} u_{i-1} - (k_{i-1} + k_{i+1}) u_i + k_{i+1} u_{i+1} \right)
\end{equation}
This allows us to write the complete discretization of (\ref{eq:integral_nd}) for spatially varying diffusion as:
\begin{align}
 h^n \sum_{j\ne i} \left[ \frac{-2(u_j - u_i) \, a_{ij}}{|\x_j-\x_i|^{n+2\beta}} + 
   \sum_{d=1}^n \overline{\partial^2_{dd} u_i} \frac{w(|\x_j - \x_i|) (\x_{j,d} - \x_{i,d})^2}{\norm{\x_j-\x_i}^{n+2\beta}} \right] \nonumber \\
  - \sum_{d=1}^n \overline{\partial^2_{dd} u_i} \int_{\Omega \cup \Omega_0} \frac{w(|\y - \x_i|) (\y_d - \x_{i,d})^2}{\norm{\y-\x_i}^{n+2\beta}} \d{\y}
\label{eq:kappa_nd}
\end{align}

Equation~(\ref{eq:kappa_nd}) reduces to (\ref{eq:I1}) and (\ref{eq:I2}) of the previous section for the one dimensional case $n=1$. The resulting multidimensional discretized operator can also be written as $A = B + D = C$ where $C$ is a sparse  symmetric matrix resulting from the treatment of the singularity. On a regular grid $C$ has a memory footprint similar to that of the discretization of the classical Laplacian with 5/7-point in 2D/3D. For the constant diffusivity case, $C$ reduces to a scaled Laplacian as derived in \cite{minden20}.

\subsection{Variable fractional order}

For the case of a spatially-varying fraction order, $\beta(y)$ can be written as the Taylor series
\begin{equation}
  \beta(\y) = \beta(\x) + \nabla \beta(\x)^T (\y - \x) + \mathcal{O}(\norm{\y - \x}^2_{\scriptscriptstyle{\nabla^2 \beta (\x)}} )
\end{equation}
where the quadratic term is the squared norm with respect to the Hessian at $\x$. This allows the singular kernel that needs to be regularized to be expressed as:
\begin{equation}
\gamma(\x, \y) = \frac{1}{\norm{\y-\x}^{n + \beta(x) + \beta(y)}} = \frac{1}{\norm{\y-\x}^{n+2\beta(\x)}} (1 - \nabla \beta(\x)^T (\y - \x) \log \norm{\y - \x} + \mathcal{O} (\norm{\y - \x)}_H^2 ) 
\label{eq:kernelexpansion_nd}
\end{equation} 
where the matrix $H$ in the quadratic term involves the Hessian $\nabla^2 \beta(\x)$ and the outer product $\nabla \beta(\x)^T \nabla \beta(\x)$. These terms do not play an explicit role in the desingularization as we explain below. 

We consider the evaluation of the integral operator (\ref{eq:integral_nd}) at a point $\x = \x_i$. Combining (\ref{eq:kernelexpansion_nd}) with 
the Taylor series expansion of $u(\y)$ of (\ref{eq:taylor_nd}), we can write 
\begin{equation}
 \mathcal{L}[u(\x_i)] = \int_{\Omega \cup \Omega_0} \left[ \frac{-2(u(\y) - u(\x_i))}{\norm{\y - \x_i}^{n + \beta(\x_i) + \beta(\y)}} + C(\y) \right] \! \d \y \,
-\int_{\Omega \cup \Omega_0}  C(y) \d \y
\label{eq:desingularized}
\end{equation}
where the desingularization term is defined as:
\begin{equation}
  C(\y) = \frac{w(\norm{\y-\x_i}) \, \left( (\y-\x_i)^T \nabla^2 u_i (\y-\x_i)  - 2 (\y-\x_i)^T \nabla \beta(\x_i) \nabla u(\x_i)^T (\y-\x_i) \log \norm{\y-\x_i} \right) }   
   {\norm{\y-\x_i}^{n+2\beta(\x_i)}}
   \label{eq:Cy}
\end{equation}
The first integral of (\ref{eq:desingularized}) is no longer singular. Its integrand is zero as $\y \to \x_i$ with enough regularity in its derivatives to admit a second-order accurate discretization by a trapezoidal rule. We point out that the higher order terms in (\ref{eq:kernelexpansion_nd}) and (\ref{eq:taylor_nd}) were not included in $C(\y)$ because they either contribute odd terms whose integrals vanish or fourth-order smooth terms that do not affect the overall resulting smoothness of the integrand. The second integral of (\ref{eq:desingularized}) can be written as the sum of two terms, each involving the product of first or second derivatives of $u$ at $\x_i$ with an integral independent of $u$ that can be carried out by a separate numerical quadrature.

As in Section~\ref{sec:varkappa_nd}, a further simplification can be performed by noting that terms involving products of the different components of the $n$-dimensional vector $(\y - \x_i)$ result in integrals that vanish. The numerator of (\ref{eq:Cy}) can be then written as the sum of $n$ terms, each involving derivatives of $u_i$ and $\beta_i$ in only one of the spatial dimensions:
\begin{equation}
  C(\y) = \sum_{d=1}^n 
  \frac{w(\norm{\y-\x_i}) \, \left( \partial^2_{dd}u_i \, (\y_d - \x_{i,d})^2 
                        - 2 \, \partial_d u_i \, \partial_d \beta_i \, (\y_d - \x_{i,d})^2 \log \norm{\y-\x_i}  \right)}
       {\norm{\y-\x_i}^{n+2\beta(\x_i)}}
\end{equation}
Plugging into (\ref{eq:desingularized}) results in the final discretization of (\ref{eq:integral_nd}) at $\x_i$:
\begin{align}
& h^n \sum_{j\ne i} \left[ \frac{-2(u_j - u_i)}{\norm{\x_j - \x_i}^{n + \beta_i + \beta_j}} \phantom{\sum_{d=1}^n} \right. \nonumber
\\ 
& + \left. \sum_{d=1}^n \partial^2_{dd} u_i \frac{w(\norm{\x_j-\x_i}) (\x_{j,d} - \x_{i,d})^2}{\norm{\x_j - \x_i}^{n + 2\beta_i}}
    -2  \sum_{d=1}^n \partial_{d} u_i \frac{w(\norm{\x_j-\x_i}) \partial_d \beta_i (\x_{j,d} - \x_{i,d})^2 \log \norm{\x_j - \x_i}}{\norm{\x_j - \x_i}^{n + 2\beta_i}}
 \right] \nonumber
 \\
& - \sum_{d=1}^n  \partial^2_{dd} u_i  \int_{\Omega \cup \Omega_0} \frac{w(\norm{\y-\x_i}) (\y_d - \x_{i,d})^2}{\norm{\y - \x_i}^{n + 2\beta_i}} \d{\y} \nonumber
 \\
& + 2 \sum_{d=1}^n  \partial_d u_i  \int_{\Omega \cup \Omega_0}  \frac{w(\norm{\y-\x_i}) \partial_d \beta_i (\y_d - \x_{i,d})^2 \log \norm{\y - \x_i}}{\norm{\y - \x_i}^{n + 2\beta_i}}  \d{\y}
\label{eq:beta_nd}
\end{align}

%

\section{Matrix representation and solution strategies} 

In this section we describe computational strategies for storing and factoring the discretized operators resulting from (\ref{eq:kappa_nd}), (\ref{eq:beta_nd}), or multi-dimensional versions of (\ref{eq:nonsym_disc}). As mentioned in Section~\ref{sec:intro}, representing these operators in their natural dense form is prohibitive because of the $O(N^2)$ memory footprint that would be required. Using the fact these matrices are ``data sparse'', i.e.,  blocks of $A$ admit low rank approximations, allows substantial reductions in memory to be realized. In Section~\ref{sec:blr}, we describe a blockwise low rank matrix representation and motivate its use, and in Section~\ref{sec:cholesky} we show that the representation allows for an efficient left-looking block Cholesky algorithm to operate directly on the compressed format and generate a factorization in $O(N^2)$ operations.

\subsection{Blockwise low rank matrix representation}
\label{sec:blr}

Consider a matrix block $A_{ts}$ where $t$ are row and column indices that correspond to clusters of grid points in spatial regions $\Omega_t$ and $\Omega_s$, respectively. $\Omega_t$ and $\Omega_s$ may conveniently be taken as axis-aligned bounding boxes of the respective point sets.  Standard results \cite{borm10} regarding the approximability of nonlocal operators and their inverses (or the closely related Cholesky factors) \cite{hackbusch15,bebendorf03,bebendorf07}, which also apply to fractional Laplacians \cite{karkulik19}, establish that if the admissibility condition:
\begin{equation}
  \max \{ \text{diam}(\Omega_t), \text{diam}(\Omega_s)\} \le \eta \, \text{dist} (\Omega_t, \Omega_s)
  \label{eq:admissibility}
\end{equation}
holds, then the asymptotically smooth kernel $\gamma(\x, \y)$ may be approximated on the bounding regions $\Omega_t$ and $\Omega_s$ by a tensor product interpolating polynomial of degree $d$ in each spatial dimension, $\bar \gamma_{ts}(\x, \y)$, with an approximation error bounded as \cite{hackbusch02}:
\begin{equation}
  \abs{\gamma(\x, \y) - \bar \gamma_{ts}(\x, \y)} \le \frac{C}{\text{dist}(\Omega_t, \Omega_s)^\sigma} q^d
  \label{eq:kernel_err}
\end{equation}
where $\sigma$ is the order of the kernel singularity, i.e., $n+2\beta$ and $n+1+\beta$ in the symmetric and non-symmetric formulations of Section~\ref{sec:formulation}, respectively, and $q = c \eta / (c \eta +1) < 1$ depends on the admissibility parameter $\eta$ which is small when $\Omega_t$ and $\Omega_s$ are well separated and grows as they become closer to each other relative to their size. $C$ and $c$ are positive constants.

A consequence of (\ref{eq:kernel_err}) is that a matrix block $A_{ts}$ of size $m \times m$ may be represented by $m \times k$ factors $U_{ts}$ and $V_{ts}$ where the block rank is $k \le d^n$, and the approximation has the error bound \cite{borm10}:
\begin{equation}
  \norm{A_{ts} - U_{ts} V_{ts}^T}_F \le \frac{C_0 \abs{\Omega_t}^{1/2} \abs{\Omega_s}^{1/2}}{\text{dist}(\Omega_t, \Omega_s)^\sigma} q^d
  \label{eq:block_error}
\end{equation}

Similar approximations bounds can also be written in the case of nested bases, i.e., when the low rank block factorization of $A_{ts}$ is expressed as $U_t S_{ts} V_s^T$ and the $U_t$ and $V_s$ bases are expressed in terms of the bases of children clusters of $t$ and $s$. 

\begin{figure}
  \begin{center}
    \includegraphics[width=.31\textwidth]{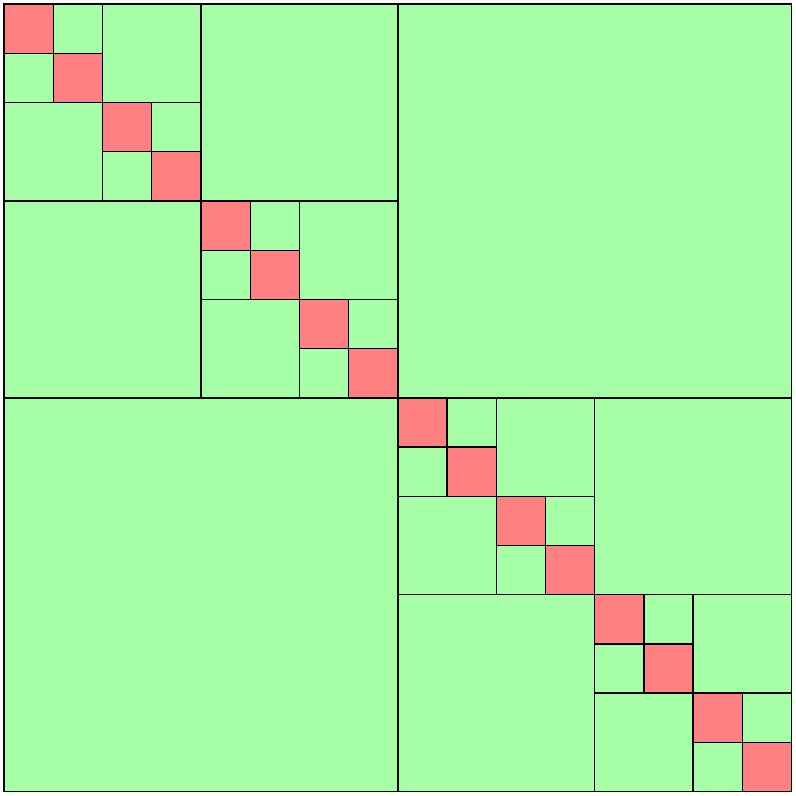}
    \hspace*{5pt}
    \includegraphics[width=.31\textwidth]{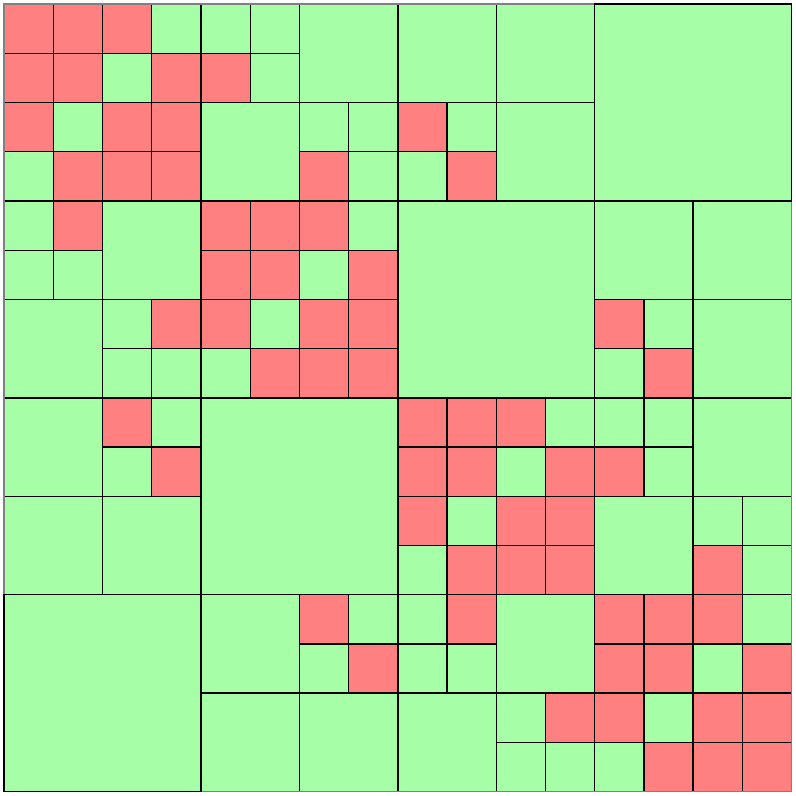}
    \hspace*{5pt}
    \includegraphics[width=.31\textwidth]{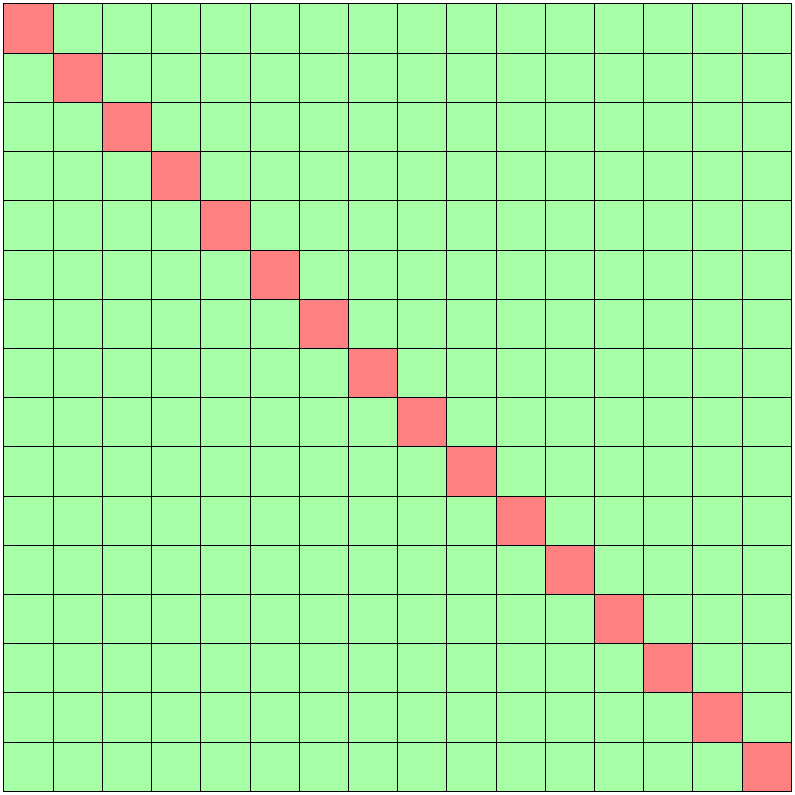}
  \end{center}
  \caption{Matrix structures that exploit local low rank approximations. Red blocks are stored as dense matrix blocks, while green blocks are stored as low rank factorizations. From left to right: weak-admissiblity $\mathcal{H}$, strong-admissibility $\mathcal{H}$, and blockwise low rank representations.}
  \label{fig:h_blr}
\end{figure}

Data sparsity may be exploited in the matrix representation in a variety of manners. A popular way of doing so has been through different flavors of hierarchical ($\mathcal{H}$) matrices. The left panel of Fig.~\ref{fig:h_blr} depicts one the simplest such representations with fixed matrix blocking, where every off-diagonal block touches the diagonal and is stored as a low rank factorization. This structure is alternatively known as  weak-admissibility $\mathcal{H}$-matrix structure, HODLR, or HSS \cite{xia10} in the case when the low rank block factorizations are expressed in nested column and row bases.  However the weak admissibility condition is only adequate for essentially one-dimensional problems and in that case direct solvers are possible \cite{gillman12}. Using the weak admissibility structure for multidimensional problems however would require very large ranks, that grow as a (fractional) power of $N$, to reach reasonable accuracy requirements since in such problems the factor $q$ that appears in (\ref{eq:block_error}) is close to 1, and the distance between point clusters is on the order of the grid spacing for the largest off-diagonal blocks.  

The middle panel of Fig.~\ref{fig:h_blr} depicts a representation that remedies this rank growth problem by allowing refinements of all blocks of the matrix, adaptively, as needed. In this representation, the dense blocks may appear anywhere in the matrix, not just along diagonal blocks, and various blocks are sized in such a way to allow bounded ranks $O(1)$ to be used everywhere. This representation is known as a standard or strong-admissibility $\mathcal{H}$-matrix, or $\mathcal{H}^2$ in the case of nested bases, and results in optimal storage complexities of $O(k N \log N)$ and $O(k N)$, respectively, with $k$ being a relatively small representative local block rank.  Unfortunately, the generality of this representation does not allow for efficient direct factorization algorithms that can be executed on multicore hardware. Therefore, general $\mathcal{H}$ and $\mathcal{H}^2$ representations have primarily been used as the main workhorse for iterative solution methods because matrix-vector multiplication can be performed efficiently \cite{boukaram19a} with them.

An alternative representation, and the one we adopt here, is a blockwise low rank representation. Instead of a full hierarchy of levels, this representation introduces only one level between the scalar operations and the full matrix dimension. Blocks are uniform in size and all off-diagonal blocks are stored as low rank factorizations, as depicted in the right panel of Fig.~\ref{fig:h_blr}. The block ranks are computed adaptively so that a uniform accuracy is maintained in all blocks. Blocks that correspond to well-separated clusters will require a small rank, while larger ranks are needed in blocks with clusters resulting in large admissibility constants $\eta$ in (\ref{eq:admissibility}). In practical 2D and 3D problems, the number of blocks requiring large ranks is quite small. Even though this representation does not attain the optimal asymptotic memory footprint of strong-admissibility $\mathcal{H}$ representations, substantial compression is achieved on problems of interest with its $O(k N^{1.5})$ asymptotic growth because of the relatively small average ranks that can be achieved, as we show in Section~\ref{sec:results}.  Additionally, by using the matrix blocks as atomic computational tiles, efficient direct factorization algorithms that benefit both from data sparsity and rich parallelism are possible, as we describe in Section~\ref{sec:cholesky} below.

In order to produce small ranks in this tile low rank (TLR) representation, a proper ordering of the grid points is essential. Ideally, points with indices close to each other should be spatially clustered together to allow (\ref{eq:admissibility}) to be satisfied with a small $\eta$. A small $\eta$ results in a small $q$ in (\ref{eq:block_error}), and therefore in smaller ranks for a given target approximation accuracy. 
A natural ordering of points in a regular grid will not, for example, satisfy this requirement. Optimal orderings for minimizing ranks are generally not known nor are practical. Instead, we use ordering heuristics similar to those developed for clustering in hierarchical matrices \cite{boukaram19a}.

We first fix the tile size $m$, which can be tuned to the cache size of the target hardware. 
The ordering of the geometric data is then determined by partitioning the grid points using a KD-tree, with repeated plane splits along coordinate directions, aimed to partition the points into clusters that are as close to the chosen tile size as possible. The construction is recursive starting from the whole point set as the topmost cluster. The points within each cluster are sorted by projecting along the largest dimension of its bounding box and then split into a left cluster whose size is half the closest power of two of the full cluster multiplied by the tile size and a right cluster containing the remaining points. This produces a cluster tree whose leaves are all the same size with the possible exception of the right most leaf, allowing the construction of the tile low rank matrix with just the final block row and column requiring padding. The resulting ordering of the grid points provides the structure and the starting point for constructing the matrix approximation and its factorization as we describe next.

\subsection{Matrix factorization in the blockwise compressed TLR format}
\label{sec:cholesky}

The first step in the processing is to construct the TLR matrix approximation. We perform this for every block/tile independently and concurrently. The $m \times m$ tiles of the matrix are evaluated, with each entry requiring a kernel evaluation and, when appropriate, a singularity correction as described in Section~\ref{sec:singularity_nd}. Each tile is then compressed using an adaptive randomized approximation (ARA) algorithm \cite{halko11,boukaram19b}. ARA requires only the sampling of the block being compressed via multiplication with random vectors. A non-adaptive randomized method generates a fixed rank $k$ approximation by: (1) sampling using the product $Y = A_{ts} X$  where $X$ is a set of $k$ random vectors, and (2) orthogonalizing $Y = QR$ to produce an approximate basis $U=Q$ for the columns of $A_{ts}$. The block is then projected onto this basis to produce the right low rank factor $V = A^T U$, thus producing the low rank factorization of $A_{ts} \approx UV^T$. Adaptive methods that automatically detect the appropriate rank for a given target accuracy, sample the matrix block $A_{ts}$ one vector at a time, iteratively constructing the orthogonal basis $U$ until the convergence threshold $\norm{A - UV^T} \le \epsilon$ is satisfied. Efficient and cache-friendly implementations of ARA are possible, and have been developed for GPU execution as well \cite{boukaram19b}.

\begin{algorithm}[H]
\caption{Left Looking Cholesky}
\label{alg:left_cholesky}
\begin{algorithmic}[1]
\Procedure{lchol}{$A, m$}
	\State $nb$ = size$(A) / m$   \Comment{\emph{number of tiles per block column}}
	\For{$k = 1 \rightarrow nb$} 
		\For{$j = 1 \rightarrow k-1$}		\label{alg:left_cholesky:lru_start}
			\For{$i = k \rightarrow nb$}
				\State $A(i, k) = A(i, k) - L(i, j) L(k, j)^T$
			\EndFor                           \label{alg:left_cholesky:lru_end}
		\EndFor
		\State $L(k, k) = $ chol$(A(k, k))$
		\For{$i = k+1 \rightarrow nb$}    
			\State $L(i, k) = A(i, k) / L(k, k)^T$
		\EndFor
	\EndFor
\EndProcedure
\end{algorithmic}
\label{alg:cholesky}
\end{algorithm}

Cholesky factorization of the constructed TLR matrix starts with a block factorization algorithm and operates on off-diagonal tiles using their low rank $U V^T $ representations. Algorithm~\ref{alg:cholesky} is a high level description of a left-looking variant of Cholesky, that operates on tiles of size $m\times m$ and updates tiles in the $k$th column using low rank updates to its left (lines 4--6). The left looking Cholesky variant has the convenient property that each tile is updated only once during execution. This is important from a performance viewpoint as it minimizes the number of tile compressions that have to be performed. Algorithms that update tiles multiple times incur additional costs due to the repeated tile compressions that would be necessary to prevent increase in ranks during intermediate computations. 

In addition, the key update operation of lines 4--6
\begin{equation}
  A(i, k)= A(i, k) - \sum_{j=1}^{k-1} L(i, j) L(k, j)^T = A(i, k) - \sum_{j=1}^{k-1} U(i, j)V(i,j)^T V(k, 
  j) U(k,j)^T 
\end{equation}
is performed using an ARA operation, i.e., its right hand side is sampled with random vectors as needed to approximate the output to the target accuracy, allowing substantial parallelism in the process.  In total, the factorization of the constructed TLR matrix can be done in $O(k N^2)$ if a block size of size $m = \sqrt{N}$ is used. 

We also note here that beyond the savings in operation count from working directly with the compressed low rank representation, the small size of the low rank data and the regularity of the tile size allow more effective use of small cache memories. The savings in latency from having the low rank data reside high on the memory hierarchy produce a significant performance boost because modern hardware architectures are provisioned for high processing power relative to memory capacity and memory bandwidth \cite{keyes20}. The numerical results in the next section show the substantial effects of these combined savings on the performance of the factorization.

\section{Numerical Results}
\label{sec:results}

Is this section we describe numerical experiments in 1D and 2D to illustrate the effectiveness of the singularity subtraction technique in the variable coefficient case and the computational savings realized by the TLR format in handling the discretized operators. The code for reproducing these experiments will be available in a branch of the H2Opus software distribution \verb|https://github.com/ecrc/h2opus|.  

\subsection{Examples in 1D}

\begin{figure}[ht]
  \begin{center}
    \includegraphics[width=.4\textwidth]{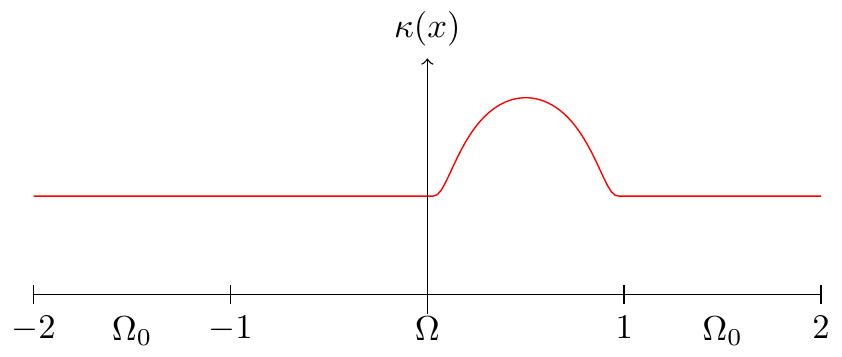}
  \end{center}
  \caption{Spatially varying nonlocal diffusion coefficient $\kappa(x)$.}
  \label{fig:profile1}
\end{figure}

We first consider the variable diffusivity case where the simulation region consists of an interior region $\Omega$, $-1 \le x \le 1$, with homogeneous ``Dirichlet'' conditions imposed outside $\Omega$ in the region $-2 \le x \le 2$. The diffusion coefficient is defined as $\kappa(x) = 1 + \mathrm{bump}(x; 0.5, 1.0)$,  where a ``bump'' function with support $\ell=1.0$ centered at $c=0.5$, and defined as: 
\begin{equation}
  \mathrm{bump}(x; c, \ell) = 
  \begin{cases}
    \mathrm{exp}\big(\! - \! \frac{1}{1-r^2}\big), \, r = \frac{x - c}{\ell /2}, & |r| < 1 \\
    0, & |r| \ge 1
  \end{cases}
\end{equation}
is added to a uniform background $\kappa = 1$. The function $\kappa(x)$ is plotted in Fig.~\ref{fig:profile1}. A constant fractional order $\beta = 0.75$ is used. 

We solve the problem for a uniform right hand side $f(x) = 1$, on regular grids of size $N$ = 64, 128, 256, 512, 1024, and 2048 in $\Omega$. Since there is no readily available analytical solution for this spatially varying coefficient case, we estimate the error on each grid by using the next finer grid as the reference solution, $e^h = u^h - u^{h/2}$, and use the sequence of error estimates to compute the rate of convergence $p$ of the discretization, $\norm{e^h} / \norm{e^{h/2}} = 2^p$. 

\begin{figure}
  \begin{center}
    \includegraphics[width=.45\textwidth]{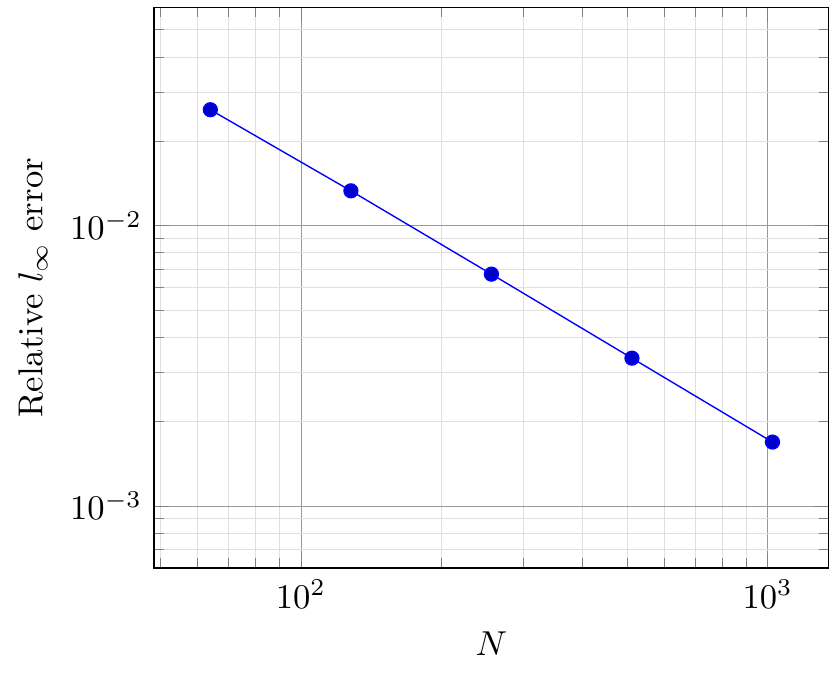}
    \includegraphics[width=.45\textwidth]{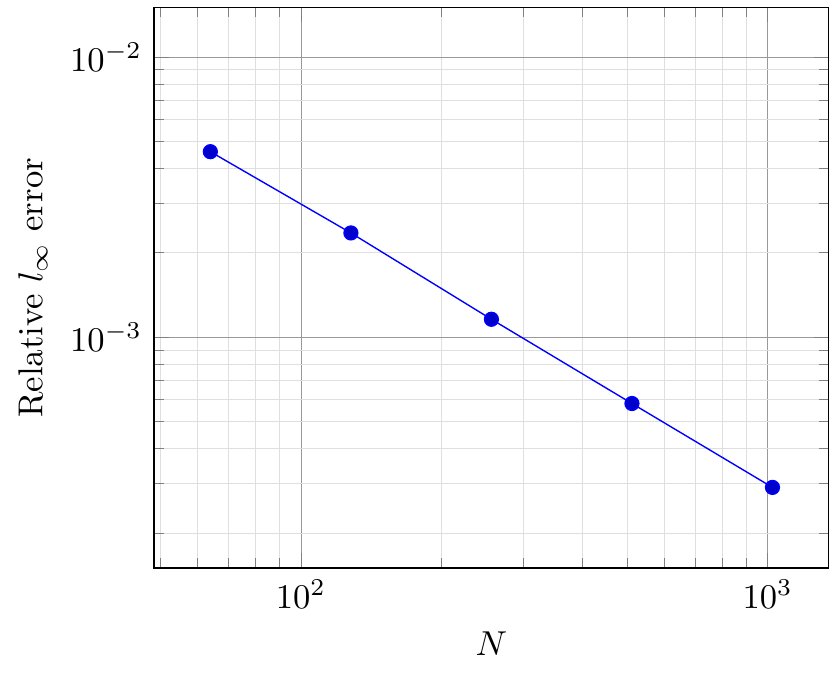}
  \end{center}
  \caption{$O(h)$ decrease in discretization error for sample problems with spatially varying coefficients in 1D: (left) variable diffusivity $\kappa(x)$, (right) variable fractional order $\beta(x)$.}
  \label{fig:var1}
\end{figure}

The left panel of Fig.~\ref{fig:var1} shows the resulting linear decrease in the relative max-norm error with grid size $h = 1/N$. We note that the first-order accuracy reached is limited only by the reduced regularity of the solution itself, which has singular derivatives at the boundaries. The trapezoidal rule with the singularity treatment can achieve second order accuracy if the solution had more regularity.


Next, we consider a variable fractional order example defined in the same interior region $\Omega$, $-1 \le x \le 1$, also with homogeneous ``Dirichlet'' conditions imposed outside $\Omega$ in the region $-2 \le x \le 2$. A linear spatial variation in $\beta$ is used, $\beta(x) = 0.7 + 0.1 x$, with a constant diffusivity coefficient $\kappa = 1$. As in the previous example,  we solve the problem for a uniform right hand side $f(x) = 1$, on grids of size $N$ = 64, 128, 256, 512, 1024, and 2048 in $\Omega$.  The relative max-norm error is computed from the difference of two solutions on successive grids. 
The right panel of Fig.~\ref{fig:var1} shows the resulting linear decrease in the error, computed as the max-norm of the difference in two solutions on successive grids.

\begin{figure}
  \begin{center}
    \includegraphics[width=.85\textwidth]{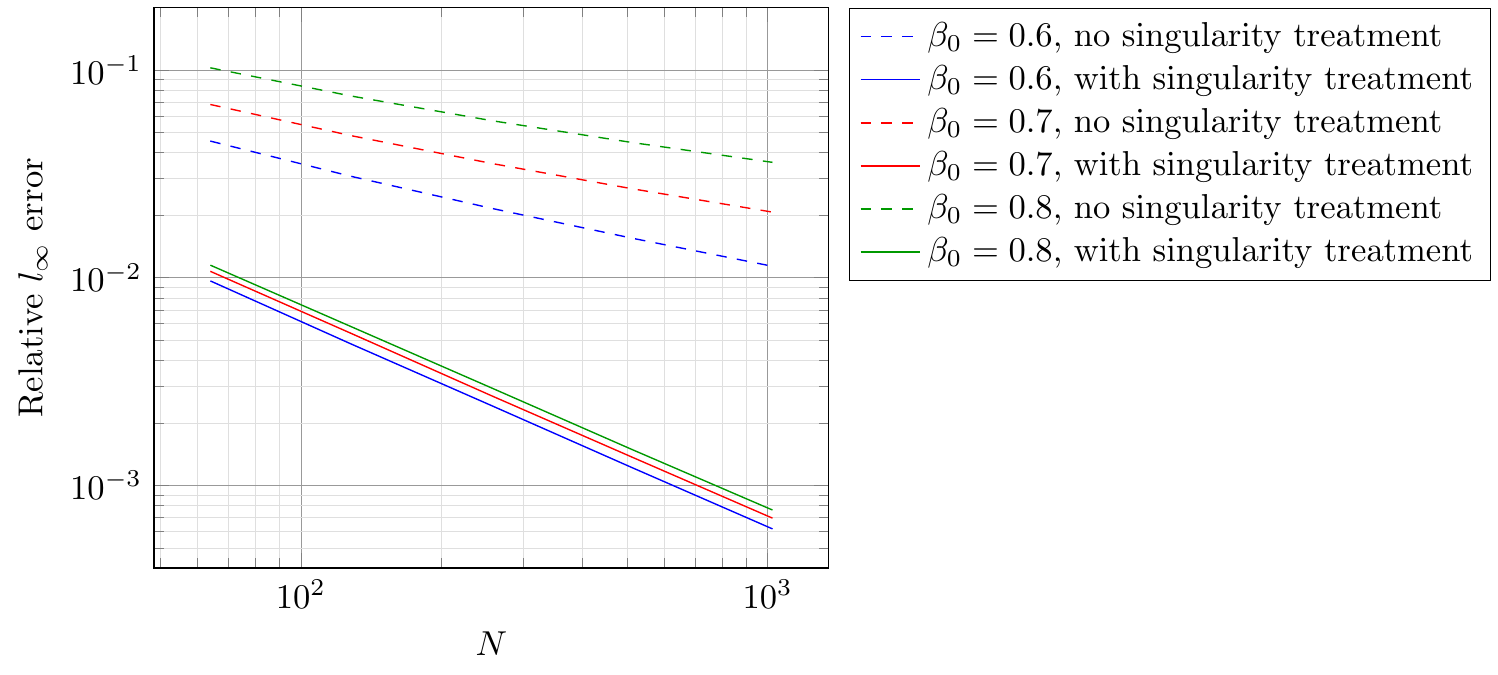}
  \end{center}
  \caption{Convergence behavior for variable fractional order of the form $\beta(x) = \beta_0 + 0.1 x$ in the non-symmetric formulation. For comparison, the behavior without the singularity treatment is also shown.}
  \label{fig:var_ns}
\end{figure}

A similar convergence behavior is obtained for the non-symmetric formulation. Figure~\ref{fig:var_ns} shows the 
convergence for problems with spatial variation $\beta(x) = \beta_0 + 0.1 x$ for three different values of $\beta_0$. In all cases, the singularity treatment results in solution convergence that is linear. For reference, the plot also shows the much slower convergence that results without the explicit treatment of the singularity. We also note that in this case the rate of convergence deteriorates faster as $\beta \to 1$, particularly as the grid is refined. 

\subsection{Variable Coefficients in 2D}
In the 2D experiments, we consider the region $\Omega = [-1, 1]^2$ extended to $[-2, 2]^2$ where homogeneous Dirichlet conditions are applied. 

\begin{figure}
  \begin{center}
    \includegraphics[height=.4\textwidth]{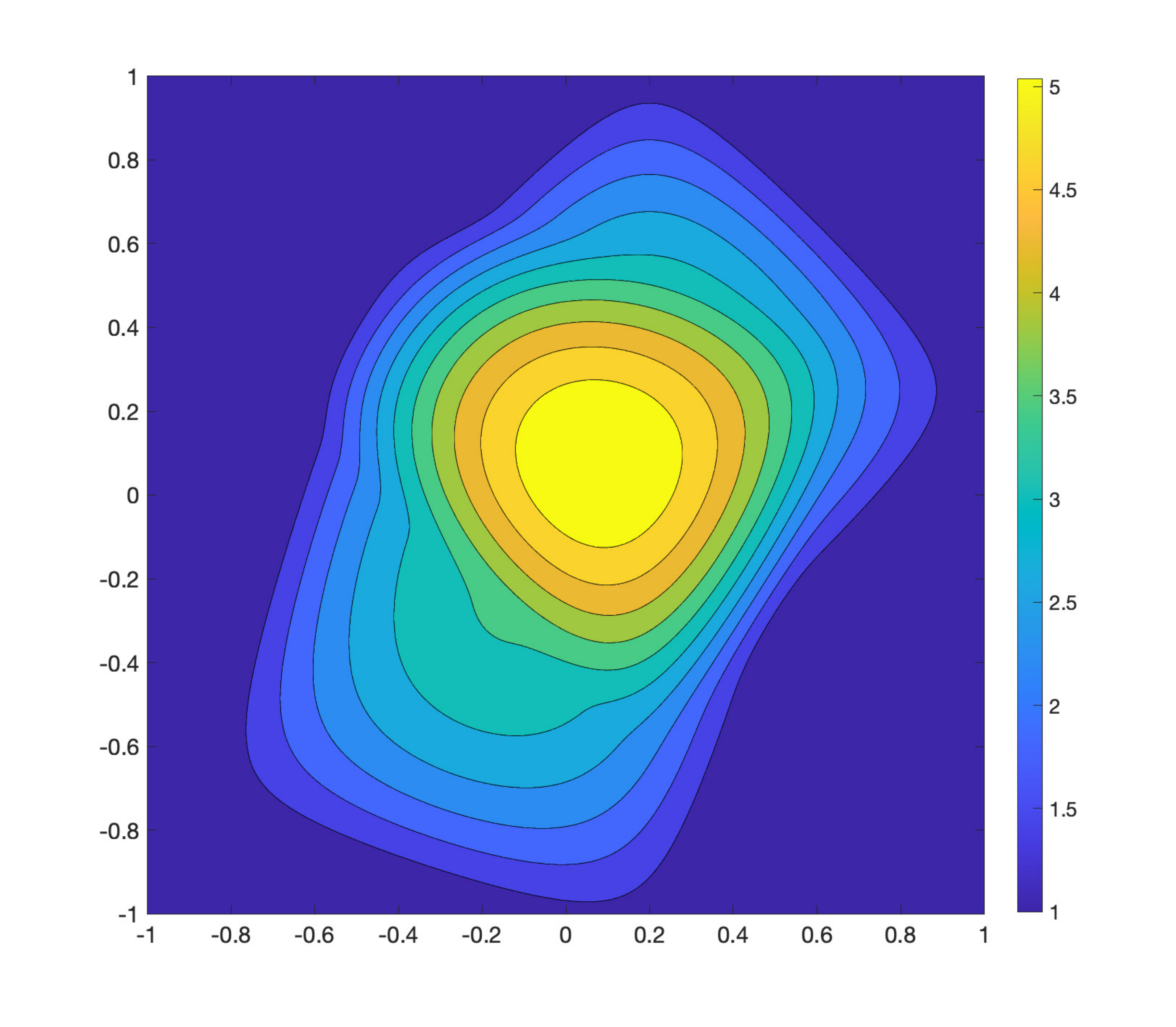}
    \includegraphics[height=.4\textwidth]{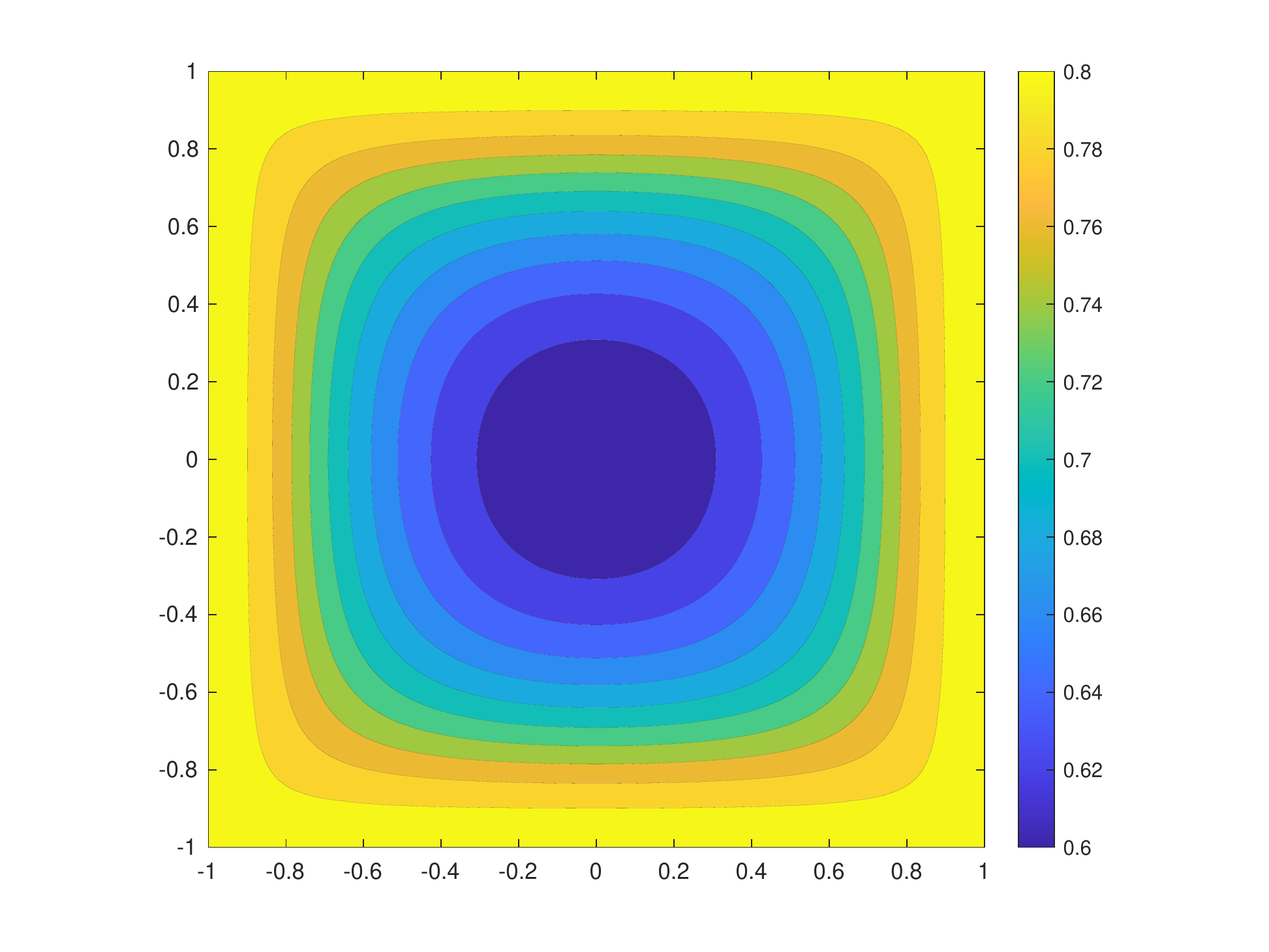}
  \end{center}
  \caption{2D spatially variable nonlocal diffusion coefficients used in the numerical experiments:
   (left) variable $\kappa(\x)$, (right) variable $\beta(\x)$.}
  \label{fig:2druns}
\end{figure}

A spatial variation of diffusivity is defined as:
\begin{equation}
  \kappa(\x) = 1 + 2.5 \, \text{bump}_{2D}(\x, \mathbf{c}_1, \bm{\ell}_1, \theta_1) + 2.5 \, \text{bump}_{2D}(\x, \mathbf{c}_2, \bm{\ell}_2, \theta_2)
\end{equation}
where $\text{bump}_{2D}(\x, \mathbf{c}, \mathbf{s}, \theta)$ is a 2D bump function obtained by taking the product of two bump functions in one variable and rotating the result by an angle $\theta$. We use $\mathbf{c}_1 = [0.2, 0.25]$, $\mathbf{c}_2 = [-0.1, -0.2]$, $\bm{\ell}_1 = [1.4, 1.4]$, $\bm{\ell}_2 = [1.4, 1.8]$, $\theta_1 = \pi/4$, and $\theta_2 = -\pi/10$. The variation is shown in the left panel of Fig.~\ref{fig:2druns}.

A spatial variation in fractional order is defined as:
\begin{equation}
  \beta(\x) = 0.8 - 0.2 \, \text{bump}_{2D}(\x, [0, 0], [2, 2], 0)
  \label{eq:beta2d}
\end{equation}
and is plotted in the right panel of Fig.~\ref{fig:2druns}. 

\begin{figure}
  \begin{center}
    \includegraphics[width=.45\textwidth]{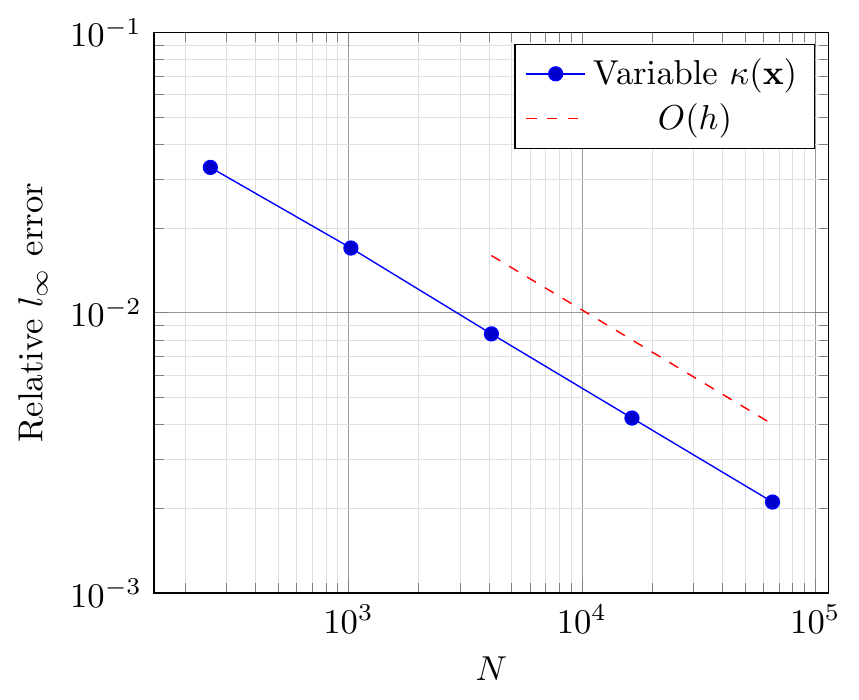}
    \includegraphics[width=.45\textwidth]{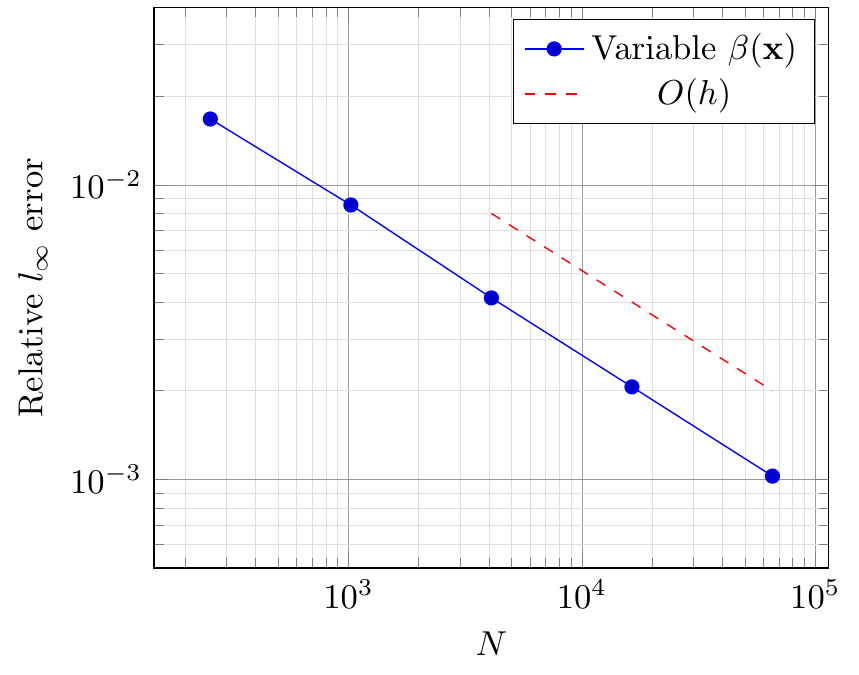}    
  \end{center}
  \caption{Convergence results for variable $\kappa(\x)$ and variable $\beta(\x)$. $O(h)$ convergence is observed in both cases, where $h$ is the grid spacing.}
  \label{fig:2dconvergence}
\end{figure}

For both variations, we solve the problem for a uniform right hand side $f(\x) = 1$ on a regular cartesian grid in $\Omega$ of size $16^2$, $32^2$ $64^2$, $128^2$, $256^2$, and $512^2$. The error is estimated on each grid by using the next finer grid as the reference solution. The relative max-norm of the error is plotted as a function of the number of grid points and the plots shown in Fig.~\ref{fig:2dconvergence}.  In both cases, we obtain the optimal rate allowed by the regularity (or rather, lack thereof) of the solution. The error decreases as $O(1/h)$ where the grid spacing is $h = N^{1/2}$.  

\begin{figure}
  \begin{center}
    \includegraphics[width=.48\textwidth]{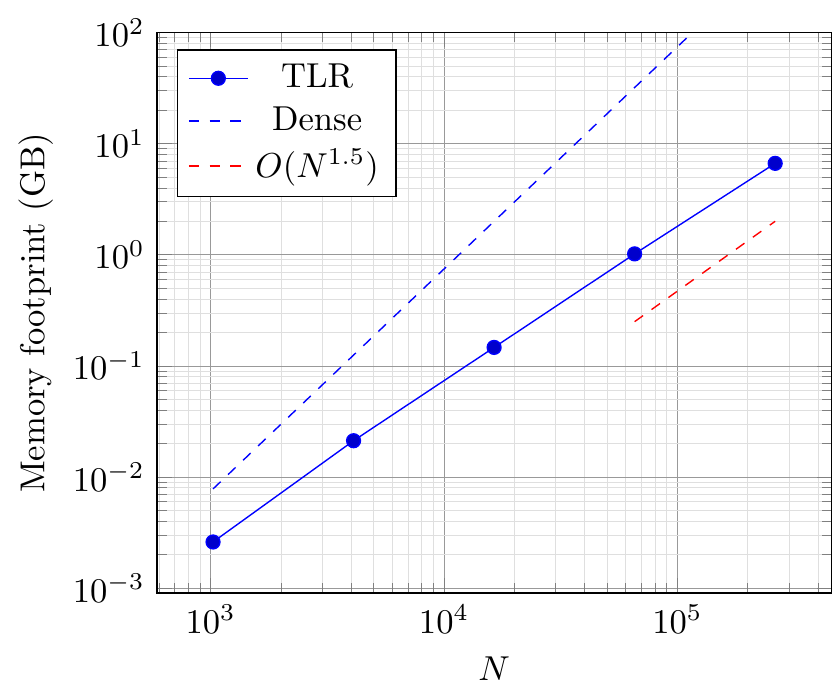}  
    \includegraphics[width=.48\textwidth]{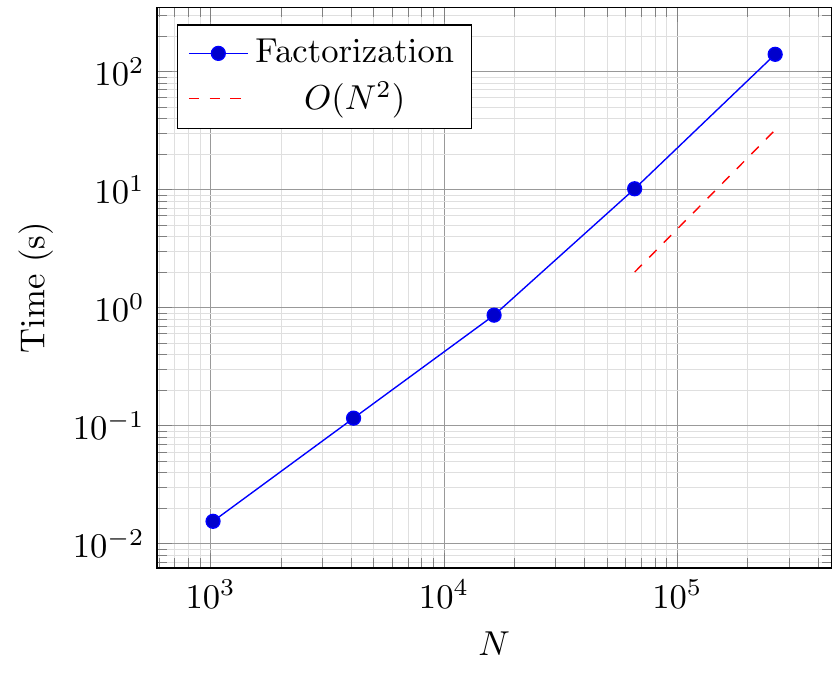}
  \end{center}
  \caption{Memory footprint and Cholesky factorization cost of the discrete operator in TLR format.}
  \label{fig:2dperf}
\end{figure}

The memory usage of the discrete operator is shown in the left panel of Fig.~\ref{fig:2dperf}. A dense representation would store $N^2$ numbers. For the simulation of size $N = 512^2 = 262K$, this would require an impractical 500+GB of storage for double precision floating point numbers. By contract, the TLR memory consumption, which approximates the matrix to an accuracy $\epsilon = 10^{-6}$ ($\norm{A_{ts} - A_{ts,\text{TLR}}} \le 10^{-6}$ for all tiles), grows at a much more modest $O(N^{1.5})$. For the simulation of size $N = 512^2$, it requires a total of only $6.65$GB of memory, with $2.0$GB for the dense diagonal tiles (256 tiles of size $m=1024$ each) and $4.65$GB for the off diagonal tiles, with an average tile rank $k$ less than 10. A KD-tree was used to decompose the grid point set recursively, with the leaves of the decomposition producing the grid ordering and the point clusters that define the matrix tiles. The tile ranks and therefore the memory consumption can be somewhat controlled by the desired target accuracy of the TLR representation, since the tile ranks are expected to change slowly with $\epsilon$, as $O(\abs{\log \epsilon}^{n+1})$.

The right panel of Fig.~\ref{fig:2dperf} is a plot of the factorization time for the Cholesky decomposition. The computations were performed in the TLR format to an accuracy of $\epsilon = 10^{-6}$ as well, on a workstation with two Xeon 20-core processors. For the $N = 512^2$ problem, the factorization required 140s. More importantly, the asymptotic growth in runtime is only $O(N^2)$, a substantial improvement over the $O(N^3)$ that would be needed for the decomposition in the dense format.  Given the triangular decomposition, a pair of forward and backward passes for computing a solution for a new right hand side only takes a small fraction of a second in the TLR format for the $N = 512^2$ problem.
  
\begin{figure}[ht]
  \begin{center}
    \includegraphics[width=.8\textwidth]{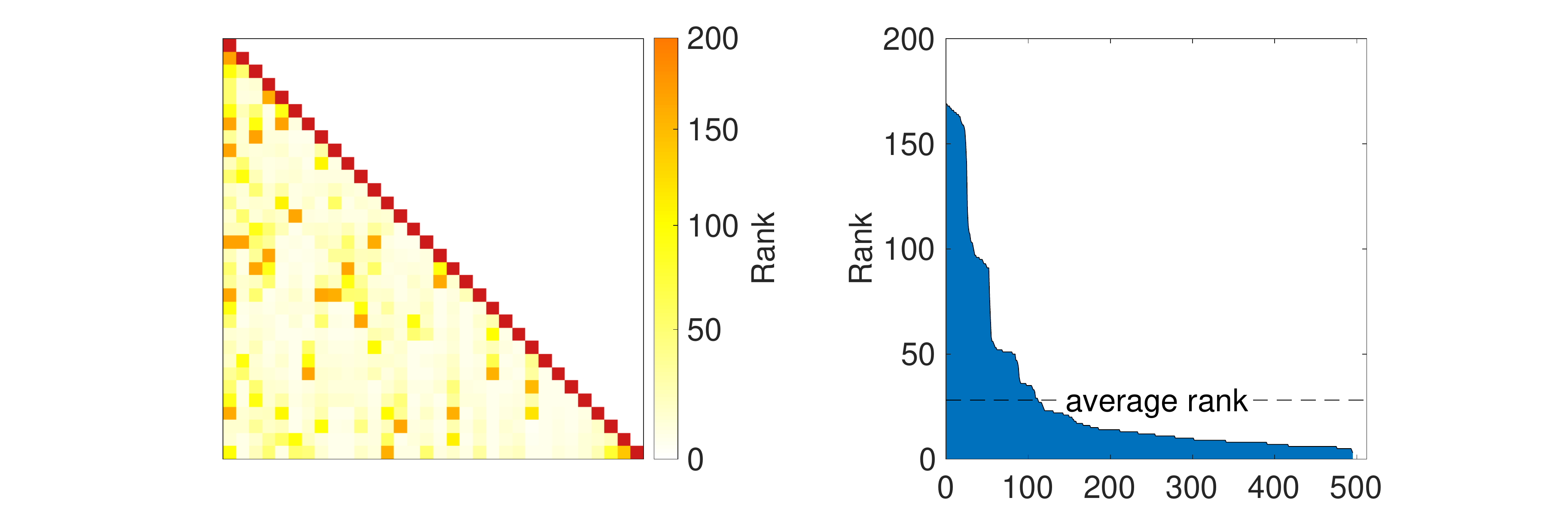}  
  \end{center}
%
  \vspace*{-25pt}
  \begin{center}
    \includegraphics[width=.8\textwidth]{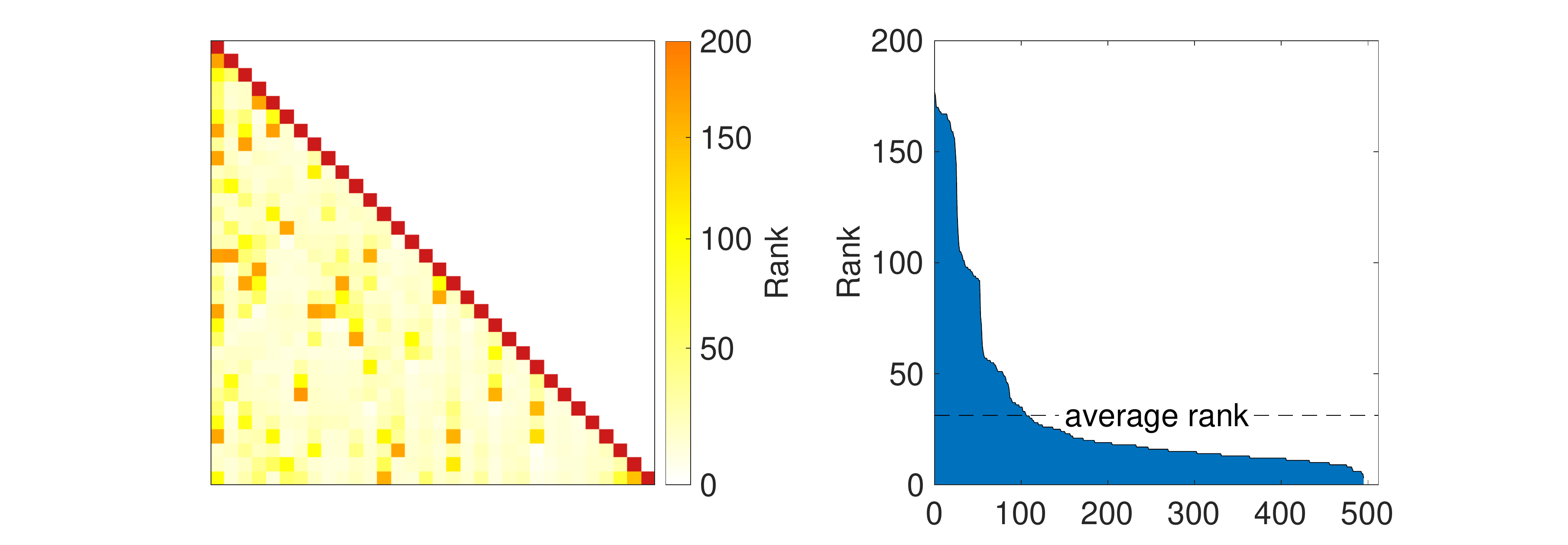}  
  \end{center}
  \caption{Tile rank distribution in the TLR format for a 16K problem with tile size $m=512$: (top) the discretized operator; and (bottom) its Cholesky factor. Left column shows heatmap plots of the matrices and right columm shows plots of their (sorted) tile ranks. Only a slight increase in tile ranks of the Cholesky factor is observed, as expected.}
  \label{fig:ranks}
\end{figure}

Finally, we plot in Fig.~\ref{fig:ranks} the rank distribution of the matrix tiles in the TLR format, both for the forward operator and its Cholesky factorization. To keep the plots legible, we use a small representative problem of size $N = 128^2 = 16K$ with the variable fractional order of (\ref{eq:beta2d}), but the trends are very similar to those of the larger problems. We use a tile size of $m=512$ which results in 32 diagonal tiles shown in red in the heatmap plots. Plots of the distribution of the ranks of the off-diagonal tiles are shown in the right column of Fig.~\ref{fig:ranks}, with an average rank of $k = 27.8$ for the forward operator. The triangular Cholesky factor increases the ranks marginally, as can be seen through the slightly darker shades of the bottom heatmap. The average rank of the off-diagonal tiles of the Cholesky factor is $k = 31.3$ and its overall memory consumption (dense diagonal plus low rank off-diagonals) increases by less than 10\% compared to the forward operator. 
We also note that the tile size $m$ provides another tuning knob to control and trade-off memory consumption vs factorization time, that is useful for high-performance contexts, but we do not discuss this fine tuning further in this work.


\section{Conclusions and future work}

We presented a singularity treatment technique that allows the effective discretization of integral formulations of variable coefficient fractional diffusion equations. A singularity subtracting term is derived at every point by matching the asymptotic singularity of the variable diffusivity and variable fractional order kernel through expansion of its various terms. The singularity is subtracted locally to produce a regular integrand that can be discretized on a regular grid by a trapezoidal rule. The resulting discretized operator is dense, however, and requires compression to make it practical for meaningful multi-dimensional simulations. To this end, we propose a tile low rank representation which partitions the dense matrix into blocks of roughly uniform size, where every off-diagonal tile is compressed and stored as its own low rank factorization. The blockwise low rank representation allows substantial compression and a much smaller memory footprint to be achieved for the fractional diffusion operators, when the grid is ordered in a way that preserves spatial proximity. 
A Cholesky decomposition algorithm operates directly on the compressed TLR representation, and uses an adaptive randomized approximation algorithm to compute the resulting tiles of the triangular factors, in a left-looking variant of the algorithm that requires only one such compression per tile. 

Numerical experiments confirm the effectiveness of the discretization. Simulations with variable diffusivity and fractional order in 1D and 2D confirm that the best convergence rate allowed by the regularity of the solution is reached.  In particular, first-order convergence is obtained for problems with singular derivatives at the boundaries. Analysis of the discretized operator in 2D problems also confirm the efficacy of the TLR representation in reducing the memory footprint from $O(N^2)$ to $O(kN^{1.5})$ with small average tile ranks when  using a KD-tree induced clustering and ordering of the grid. A nearly two-order of magnitude reduction in memory compared to a dense format is obtained for an $N = 262K$ problem compressed to a $10^{-6}$ accuracy. Results also confirm that the direct factorization of the operator can be done in $O(k N^2)$.

These encouraging results point to a number of extensions that we intend to consider in the future. We have dealt with isotropic coefficients and plan to extend the treatment to the practically important anisotropic case. We also plan to analyze the effect of discontinuous coefficients. Our discretization has been on a regular cartesian grid, but can be extended to general geometries and triangular meshes. In addition, we plan to explore the role of GPUs in accelerating the arithmetically intensive TLR computations, which we expect should give the computations a substantial performance boost. Finally, we intend to tackle large scale 3D problems which will likely require distributed-memory computers and perhaps the use of hierarchical matrices, with their optimal computational complexities at scale.


\appendix
\section{Treatment of kernel singularity in the non-symmetric formulation}
\label{app:nonsym}

In 1D the fractional flux is written as: 
\begin{align}
  Q^\beta(x) &= -\kappa(x) \frac{d^{\beta(x)}}{dx^\beta} u(x) \\
       &= -\kappa(x) \omega(x) \int_\Omega \frac{y-x}{|y-x|^{\beta(x)+2}} u(y) dy
 \label{eq:nonsym_integral}
\end{align}

Consider the desingularization of the integral of (\ref{eq:nonsym_integral}),
\begin{align}
  I(x) = & \int_\Omega \left[ \frac{y-x}{|y-x|^{\beta(x)+2}} u(y) - \frac{(y-x) \, w(y-x) \, [u(x) + u'(x) (y-x)]}{|y-x|^{\beta(x)+2}} \right] dy \\
    & + u'(x) \int_\Omega  \frac{(y-x) \, w(y-x) \, [u(x) + u'(x) (y-x)]}{|y-x|^{\beta(x)+2}} dy 
\end{align}

The first integrand is now no longer singular and can be readily discretized by a trapezoidal/midpoint rule.

Because the integrals of odd powers of $(y-x)$ evaluate to zero, the expression for $I(x)$ simplifies to: 
\begin{align}
  I(x) = \int_\Omega \left[ \frac{y-x}{|y-x|^{\beta(x)+2}} u(y) - u'(x) \frac{w(y-x)}{|y-x|^{\beta(x)}} \right] dy 
  + u'(x) \int_\Omega  \frac{w(y-x)}{|y-x|^{\beta(x)}} dy 
  \label{eq:correction}
\end{align}

Let the local function $w(y-x_i)$ have support in a small region $\delta \le y-x_i \le \delta$, where $\delta$ corresponds to a few cell widths $m h$. Then the second term of (\ref{eq:correction}) may be discretized as:
\begin{equation}
  \int_\Omega  - u'(x) \frac{w(y-x)}{|y-x|^{\beta(x)}} dy  = u'(x) \, h \sum_{j=i-m}^{i+m-1} -\frac{w((x_{j+1/2} - x_i))}{|x_{j+1/2} - x_i|^{\beta(x_i)}} =  u'(x) C_1(x_i)
\end{equation}

If we denote the last integral of (\ref{eq:correction}) by $C_2(x)$, then
\begin{equation}
  C_2(x_i) =  \int_{x_i-\delta}^{x_i+\delta}  \frac{w(y-x_i)}{|y-x_i|^{\beta(x_i)}} dy 
\end{equation}
and the final discretization of the flux becomes: 
\begin{align}
  Q^\beta(x_i) \approx -\kappa(x_i) \omega(x_i) \, \left( h \sum_j \frac{ (x_{j+1/2} - x_i) u_{j+1/2}}{|x_{j+1/2} - x_i|^{\beta(x_i) +2}}  
    +  \frac{u_{i+1/2} - u_{i-1/2}}{h}  (C_1(x_i) + C_2(x_i))  \right)
    \label{eq:discretization}
\end{align}

$C_1$ and $C_2$ can be computed first for all $x_i$ and then used in the discretization of (\ref{eq:discretization}).

The fractional diffusion operator is then
\begin{align}
  \frac{d Q^{\beta}(x_i)}{dx} & \approx \frac{1}{h} \left( Q^{\beta} (x_{i+1}) - Q^{\beta} (x_i) \right)  \\
  & = -\kappa(x_i) \omega(x_i) \left( 
    \sum_j \frac{ (x_{j+1/2} - x_{i+1}) u_{j+1/2}}{|x_{j+1/2} - x_{i+1}|^{\beta(x_{i+1}) +2}}   
     - \sum_j \frac{ (x_{j+1/2} - x_i)  u_{j+1/2}}{|x_{j+1/2} - x_i|^{\beta(x_i) +2}}   \right) \nonumber \\
     & \quad -\kappa(x_i) \omega(x_i) \frac{1}{h^2} \left(C_{i+1} u_{i+3/2}  - (C_i + C_{i+1}) u_{i+1/2} + C_i u_{i-1/2} \right)
     \label{eq:nonsym_disc}
\end{align}
where $C_i = C_1(x_i) + C_2(x_i)$.

\section*{Compliance with Ethical Standars}
\textbf{Funding:} The authors acknowledge the support of the Extreme Computing Research Center at KAUST.\\ 
\textbf{Conflict of Interest:} On behalf of all authors, the corresponding author states that there is no conflict of interest.


\bibliographystyle{spmpsci}      
\bibliography{vcfd}

\end{document}